\documentclass[a4paper]{article}
\usepackage{amsthm,amsmath,amssymb}
\newtheorem{teor}{Theorem}[section]
\newtheorem{defin}[teor]{Definition}
\newtheorem{lemm}[teor]{Lemma}
\newtheorem{osse}[teor]{Remark}
\newtheorem{prop}[teor]{Proposition}
\newtheorem{defi}[teor]{Definition}
\newtheorem{coro}[teor]{Corollary}
\newtheorem{prob}[teor]{Problem}

\newcommand{\bele}{\begin{lemm}\begin{sl}}
\newcommand{\enle}{\end{sl}\end{lemm}}
\newcommand{\bedef}{\begin{defi}\begin{sl}}
\newcommand{\eddef}{\end{sl}\end{defi}}
\newcommand{\bete}{\begin{teor}\begin{sl}}
\newcommand{\ente}{\end{sl}\end{teor}}
\newcommand{\beos}{\begin{osse}\begin{rm}}
\newcommand{\eddos}{\end{rm}\end{osse}}
\newcommand{\bepr}{\begin{prop}\begin{sl}}
\newcommand{\empr}{\end{sl}\end{prop}}
\newcommand{\bepro}{\begin{prob}\begin{rm}}
\newcommand{\empro}{\end{rm}\end{prob}}
\newcommand{\bede}{\begin{defin}\begin{sl}}
\newcommand{\edde}{\end{sl}\end{defin}}
\newcommand{\beco}{\begin{coro}\begin{sl}}
\newcommand{\enco}{\end{sl}\end{coro}}

\newcommand{\thspace}{\hspace{3mm}}

\newcommand{\qquand}{\qquad\text{and}\qquad}
\newcommand{\quext}{\quad\text}
\newcommand{\qquext}{\qquad\text}

\newcommand{\RR}{\mathbb{R}}
\newcommand{\EE}{\mathbb{E}}
\newcommand{\NN}{\mathbb{N}}

\newcommand{\beeq}[1]{\begin{equation}\label{#1}}
\newcommand{\eddeq}{\end{equation}}

\newcommand{\beeqa}[1]{\begin{eqnarray}\label{#1}}
\newcommand{\eddeqa}{\end{eqnarray}}

\newcommand{\beal}[1]{\begin{align}\label{#1}}
\newcommand{\eddal}{\end{align}}

\newcommand{\bespl}[1]{\begin{split}\label{#1}}
\newcommand{\edspl}{\end{split}}

\newcommand{\bega}[1]{\begin{gather}\label{#1}}
\newcommand{\edga}{\end{gather}}

\newcommand{\beeqax}{\begin{eqnarray*}}
\newcommand{\eddeqax}{\end{eqnarray*}}

\def\qed{\ifmmode 
  \else \leavevmode\unskip\penalty9999 \hbox{}\nobreak\hfill
  \fi
  \quad\hbox{\hskip.5em\vrule width.4em height.6em depth.05em\hskip.1em}}
\def\endproofsym{\qed}
\newcommand{\dimbox}{\hbox{\hskip.5em\vrule width.4em height.6em depth.05em\hskip.1em}}
\renewenvironment{proof}[1][Proof]{\trivlist\item[\hskip\labelsep{\hskip0pt
    {\normalfont\scshape#1.}\hskip .321429\parindent}]\ignorespaces}
{\endproofsym\endtrivlist}
\def\endnobox{\def\endproofsym{}\end{proof}\def\endproofsym{\qed}}

\newcommand{\no}{\nonumber}

\newcommand{\beeqao}{\begin{eqnarray}\no}
\newcommand{\bealo}{\begin{align}\no}
\newcommand{\besplo}{\begin{split}\no}
\newcommand{\begao}{\begin{gather}\no}

\newcommand{\duav}[1]{\langle{#1}\rangle}

\newcommand{\ppp}{{\mathfrak p}}

\newcommand{\Dt}{\partial_t}
\newcommand{\Nx}{\nabla}
\newcommand{\Dx}{\Delta}

\newcommand{\perogni}{\forall\,}

\newcommand{\itt}{\int_0^t}

\newcommand{\io}{\int_\Omega}

\newcommand{\iTT}{\int_0^T}

\newcommand{\iTo}{\iTT\!\io}

\newcommand{\epsi}{\varepsilon}
\newcommand{\ee}{_{\varepsilon}}
\newcommand{\een}{_{\varepsilon,\nu}}
\newcommand{\lla}{_{\lambda}}

\newcommand{\OO}{_{\Omega}}

\def\eb{\varepsilon}
\def\R{\mathbb R}

\newcommand{\bn}{\boldsymbol{n}}

\newcommand{\dn}{\partial_{\bn}}

\newcommand{\lhs}{left hand side}
\newcommand{\rhs}{right hand side}

\DeclareMathOperator{\tr}{tr}
\DeclareMathOperator{\dive}{div}
\DeclareMathOperator{\deriv}{d}

\DeclareMathOperator{\Id}{Id}

\DeclareMathOperator{\loc}{loc}

\newcommand{\HUH}{H^1(0,T;H)}
\newcommand{\HUV}{H^1(0,T;V)}

\newcommand{\LDH}{L^2(0,T;H)}
\newcommand{\LDV}{L^2(0,T;V)}
\newcommand{\LDVp}{L^2(0,T;V')}

\newcommand{\LIV}{L^\infty(0,T;V)}

\newcommand{\LDHD}{L^2(0,T;H^2(\Omega))}

\let\TeXchi\chi
\def\chi{{\setbox0 \hbox{\mathsurround0pt
$\TeXchi$}\hbox{\raise\dp0 \copy0 }}}

\newcommand{\calX}{{\mathcal X}}

\newcommand{\calA}{{\mathcal A}}

\newcommand{\calE}{{\mathcal E}}

\newcommand{\calB}{{\mathcal B}}

\newcommand{\calZ}{{\mathcal Z}}

\newcommand{\barO}{\overline{\Omega}}

\newcommand{\barn}{\overline{n}}

\newcommand{\dit}{\deriv\!t}

\newcommand{\ddt}{\frac{\deriv\!{}}{\dit}}

\newcommand{\eet}{_{\varepsilon,t}}
\newcommand{\ditau}{\deriv\!\tau}
\newcommand{\zee}{_{0,\varepsilon}}

\numberwithin{equation}{section}
\begin{document}

\title{Existence of solutions and separation from singularities
for a class of fourth order degenerate parabolic equations}

\author{Giulio Schimperna\\
Dipartimento di Matematica, Universit\`a di Pavia,\\
Via Ferrata~1, I-27100 Pavia, Italy\\
E-mail: {\tt giusch04@unipv.it}\\
\and
Sergey Zelik\\
Department of Mathematics, University of Surrey,\\
Guildford, GU2 7XH, United Kingdom\\
E-mail: {\tt S.Zelik@surrey.ac.uk}
}


\maketitle
\begin{abstract}
 A nonlinear parabolic equation of the fourth order is analyzed.
 The equation is characterized by a mobility coefficient
 that degenerates at $0$.
 Existence of at least one weak solution is proved by using
 a regularization procedure and deducing
 suitable a-priori estimates. If a viscosity
 term is added and additional conditions on the nonlinear terms
 are assumed, then it is proved that any
 weak solution becomes instantaneously strictly positive. 
 This in particular
 implies uniqueness for strictly positive times
 and further time-regularization properties.
 The long-time behavior of the problem is also investigated
 and the existence of trajectory attractors and, under
 more restrictive conditions, of strong global attractors
 is shown.
\end{abstract}


\vspace{2mm}

\noindent {\bf AMS (MOS) subject clas\-si\-fi\-ca\-tion:}
35K35, 35K65, 37L30.
\thspace

\vspace{2mm}

%
%
%
%


\section{Introduction}
\label{secintro}

This paper is devoted to the analysis of the following
class of fourth order parabolic equations:
\begin{align}\label{tf1}
  & u_t-\dive (b(u)\nabla w)=0,\\
 \label{tf2}
  & w = \delta u_t - \Delta u + f(u) + \gamma(u) - g,
\end{align}
on $\Omega\times (0,+\infty)$,
$\Omega$ being a bounded smooth subset of
$\RR^d$, $d\in\{2,3\}$, coupled with
the initial and boundary conditions
\begin{align}\label{iniz-intro}
  & u|_{t=0} = u_0,
  \quext{in }\,\Omega,\\
 \label{neum-intro}
  & \dn u=b(u)\nabla w\cdot \bn=0,
  \quext{on }\,\partial\Omega.
\end{align}
The function $b(u)$ represents a solution-dependent
mobility coefficient that possibly {\sl degenerates}
at $0$ as a power of $u$
(cf.~\eqref{hpb} below), while
the sum $f+\gamma$
stands for the derivative of a configuration
potential $W$. In particular, we assume that
$f$ is the (dominating) monotone part,
with $f(u)\sim u^{-\kappa}$
for some $\kappa > 1$,
and that $\gamma$ is a bounded and
globally summable perturbation that
accounts for possible nonconvexity of $W$.
The coefficient $\delta$ in \eqref{tf2}
is assumed to be nonnegative, with $\delta>0$
describing the presence of {\sl viscosity}\/
effects. Finally, $g$ is a smooth
external forcing term.

In the two-dimensional case, problem
\eqref{tf1}-\eqref{neum-intro}
can describe the evolution
of some classes of thin liquid films,
with $u$ representing the height of the film.
Then, the singular behavior of $f$ near
$0$ accounts for the presence of short-range
repulsive forces, while the nonmonotone character
of $W'$ at $\infty$ (given by the term $\gamma$)
is related to the occurrence
of long-range repulsive forces.
An extensive presentation of the
underlying physical situation is given
in \cite{BGW} to which we refer
for more details (see also \cite{BeGr,GrR2} and
Remark~\ref{ossef} below).

In the three-dimensional case, the model is
also physically relevant since it is closely
related with the Cahn-Hilliard equation
\cite{CH} with nonconstant mobility analyzed
in a number of recent papers
both in the nondegenerate and in the degenerate case
(cf., e.g., \cite{BB,BBG,EG,NCSh,S} and the references
therein). Indeed, if $\pm1$ represent the pure states
of the order parameter in the Cahn-Hilliard
model, then we can modify $f$ into the form
$f(u)\sim (1-u^2)^{-\kappa}$ and accordingly
suppose that $b$ degenerates near $\pm1$
(instead that near $0$) as a power of $(1-u^2)$.
Then, the so-modified system \eqref{tf1}-\eqref{tf2}
turns out to represent a variant of the
Cahn-Hilliard equation
with degenerate mobility and singular potentials
analyzed in the celebrated paper \cite{EG}.
Correspondingly, all the results proved
here for \eqref{tf1}-\eqref{tf2}
also apply to the Cahn-Hilliard setting
with straighforward modifications
in the proofs (to be more precise, in the
Cahn-Hilliard setting we would even have
slightly stronger results since it would no longer
be necessary to take care of the growth of
$b$ at infinity). This also motivates
the choice of considering also the viscous case
$\delta>0$, which is particularly meaningful
in the context of Cahn-Hilliard models
(cf.~\cite{Gu}, see also \cite{NC2}).

Initial-boundary value problems related to
\eqref{tf1}-\eqref{tf2} have been addressed in
a number of recent contributions. In particular,
in \cite{BGW} further qualitative properties of the
solutions are proved in the one-dimensional case
and the stability properties of the steady states
are investigated. The papers \cite{GrR2} (devoted to
the one-dimensional case) and \cite{Gr}
(considering space dimensions 2 and 3) analyze
problem \eqref{tf1}-\eqref{neum-intro} under
assumptions on the nonlinear terms very similar to ours.
In \cite{Gr,GrR2}, existence of a solution is proved
by means of a nonnegativity-preserving finite-element
scheme, which is also effective for a numerical
investigation for the model. Finally, we mention the
recent work \cite{WZ}, where the  long-time behavior
of the problem is studied in the one-dimensional
setting. In this situation, the authors
can prove strict positivity of the
solution also in the nonviscous
case, which allows them to show existence of a smooth
global attractor by relying on the standard theory
of infinite-dimensional dynamical systems.

Our first purpose in this paper is to prove existence
of at least one weak solution to
Problem~\eqref{tf1}-\eqref{neum-intro}
under general assumptions on the data,
by using a regularization -- a priori estimate --
passage to the limit procedure. Compared to the
proof given in \cite{Gr}, our method
has the advantage to be relatively simple.
Moreover, as a byproduct of our procedure 
we see that the possibly singular solutions 
to~\eqref{tf1}-\eqref{neum-intro} can be 
approximated by the smooth and {\sl positive}\/
solutions of a regular PDE. Actually, if we 
consider, for instance, the nonviscous equation 
(i.e., the case $\delta=0$),
then we can construct the $\epsi$-approximation
taking $\delta_\epsi>0$ 
(which, as noted above, is physically motivated
at least in the Cahn-Hilliard setting)
and choosing a sufficiently singular function $f_\epsi$,
with $\delta_\epsi \searrow 0$ and $f_\epsi\to f$ in 
the limit $\epsi\searrow 0$. Then, for $\epsi>0$ 
Theorem~\ref{teosep} applies (cf.~also
Remark~\ref{arbireg}); hence, the approximating
solutions $u\ee$ are smooth and positive. Moreover, the very
same argument used to pass to the limit in 
Subsec.~\ref{pass-lim} below shows that any limit point
of $\{u\ee\}$ for $\epsi\searrow0$
solves the original nonviscous problem.

In comparison with \cite{Gr,GrR2}, we have here the extra assumption
that the singularity of $1/b(u)$ at $u=0$ is not stronger than 
the singularity of the potential $F$ (the antiderivative
of $f$), cf.~\eqref{hpf} below. 
On the one hand, this assumption looks natural and 
is satisfied for the physically relevant examples of $b(u)$ and $f(u)$. 
Indeed, the most used non-linearity $b$ in the theory of thin films is
(cf., e.g., \cite{BGW} or \cite{Gr})
\begin{equation}\label{whichb}
  b(u)=u^3+\beta^{3-n}u^n,
   \qquad n\in(0,3),
\end{equation}
and for the function $f$ we have the following model examples
(cf.~\cite{BGW}):
\begin{equation}\label{whichf}
  f(u)= Au^{-3}-Bu^{-9} \quext{or }\,f(u)=Au^{-3}-Bu^{-4},
\end{equation}
where $A,B>0$. Thus, 
in all the physically relevant cases \eqref{hpf} is satisfied.
On the other hand, this assumption allows us not only to simplify
the proof of the existence of a weak solution and to consider 
more general functions $f$ and $b$, but also to use the Moser-type 
iteration technique for improving the regularity of the constructed 
weak solution in the viscous case.

Once the existence of a solution is achieved, we
show dissipativity of the (multivalued)
dynamical process associated to the system.
This in particular entails existence
of a (weak) trajectory attractor in the sense
of Chepyzhov and Vishik \cite{CV}. 

Our subsequent results regard only the viscous
case $\delta>0$ (and in particular can be
applied to the Cahn-Hilliard model up to
the modifications described above).
Then, we can also prove that the 
trajectory attractor can be intended in the
strong sense (i.e., w.r.t.~the strong topology
of the natural phase space) at least
under slightly more restrictive conditions
on~$f$. The proof relies on an ad-hoc integration
by parts formula and a variant of the so-called
energy method (cf.~\cite{BVbook}, see also
\cite{Ba1,MRW} and \cite{CVZ1,CVZ2,EKZ} 
for applications to trajectory attractors).
If we furtherly restrict the class of 
admissible functions $f$ (namely, asking
$\kappa$ to be large enough), then we 
can also prove that weak solutions 
become uniformly separated from $0$ for
any time $t>0$, so that the degenerate
character of the system is actually lost
for strictly positive times.
This result, which is in our opinion
the main achievement of this paper, is shown by
means of a suitable version of the Moser iteration scheme
which takes time regularization effects into account.
As a further consequence of this ``separation'' property,
we can also prove arbitrarily high regularity
of weak solutions, as well as uniqueness,
at least for $t>0$. In turn, this permits to interpret
the global attractor in the frame of the standard
(single-valued) theory \cite{BVbook,Te}, rather
than in the trajectory sense.

The plan of the paper is as follows. 
In the next Section~\ref{main} we will report our notation 
and hypotheses and the statement of our existence result.
Its proof is divided into several steps and will be presented
in Section~\ref{proof-esi}.
In Section~\ref{sec:weakattr}, we will show the
existence of weak trajectory attractors.
In~Section~\ref{sec:strongattra}, we will
prove existence of a strong trajectory attractor
in the viscous case $\delta>0$
by applying the so-called energy method. 
Finally, in Section~\ref{sec-long}, 
we will prove the strict positivity of
$u$ in the viscous case under more restrictive
growth conditions on $f$.


\section{Existence result}
\label{main}

Let $\Omega$ be a smooth bounded domain of $\RR^d$,
$d\in\{2,3\}$. Let $T>0$ a given final time and
let $Q:=\Omega\times(0,T)$. Let $H:=L^2(\Omega)$,
endowed with the standard scalar product $(\cdot,\cdot)$
and norm $\| \cdot \|$. Let also $V:=H^1(\Omega)$.
We note by $\|\cdot\|_X$
the norm in the generic Banach space $X$.

We make the following assumptions on data:
\begin{align}\label{hpb}
  & b(r) = r^s + \beta r^n, \qquad \beta\ge 0,~~r\ge 0;\\
  \label{hpb2}
  & 0 \le n \le s < 10~~\text{if }\,d=3, \qquad
   0 \le n \le s~~\text{if }\,d=2;\\
 \label{hpf}
  & f(r) = -\frac1{r^\kappa}, \qquad \kappa>1,~~\kappa \ge s+1,~~r > 0;\\
 \label{hpgamma}
  & \gamma\in W^{1,\infty}(\RR)\cap L^1(\RR);\\
 \label{hpg}
  & g \in V\cap L^\infty(\Omega), \qquad
   \hat g:=\|g\|_{L^\infty(\Omega)};\\
 \label{hpdelta}
  & \delta\ge 0.
\end{align}
\beos\label{sullaregdati}
 Our key assumption here is that the singular character of $f$
 has to dominate over the degeneracy of $b$ 
 at 0 (cf.~\eqref{hpf}). Actually, one could see
 with straighforward modifications in the proofs that
 in the case $\beta>0$ (i.e., if $b(r)$ has a lower order 
 degeneration $r^n$ at $0$), then it would be enough
 to ask $\kappa\ge n+1$ rather than $\kappa\ge s+1$.
 We also point out that the requirement $s<10$ in the 
 three-dimensional case is motivated by the growth
 of $b$ at $\infty$ (and not by its degeneration at $0$).
 Consequently, in the application to the Cahn-Hilliard
 model (where solutions have to stay in between
 the two barriers $r=\pm1$), $s$ could in fact 
 be arbitrary also in the case $d=3$.
\eddos
\noindent%
We also define, whenever they make sense, the following
functions:
\begin{equation}\label{defiFW}
  F(r):=\frac1{\kappa-1} + \int_1^r f(\tau)\,\ditau, \qquad
   \Gamma(r):=\int_1^r \gamma(\tau)\,\ditau, \qquad
   W(r):= F(r)+\Gamma(r).
\end{equation}
\beos\label{ossef}
 According, e.g., to \cite{BGW}, a physically relevant
 expression for $W'=f+\gamma$ is given by
 \begin{equation}\label{physf}
   W'(r) \sim  -\frac1{r^\kappa} + \frac1{r^k},
    \qquext{where }\,k\in(1,\kappa).
 \end{equation}
 Actually, this situation is not covered by our assumptions
 \eqref{hpf}-\eqref{hpgamma}. However,
 it is clear that, just with technical modifications
 in the proofs, one could replace \eqref{hpf} with
 something like
 \begin{equation}\label{hpfext}
   c_1 \frac1{r^{\kappa+1}}
    \le f'(r)
    \le c_2 \frac1{r^{\kappa+1}},
 \end{equation}
 for all $r>0$ and some $c_1,c_2>0$.
 Assuming \eqref{hpfext} and properly
 choosing $\gamma$, it is clear that we can
 deal with the case \eqref{physf}.
 Nonetheless, we will assume \eqref{hpf}
 in place of \eqref{hpfext} in order to
 reduce technical complications in the proofs.
\eddos
Next, we
%
%
introduce the energy functional associated to
system \eqref{tf1}-\eqref{tf2}:
\begin{equation}\label{defiE}
  \calE(u):=\io\Big(
   \frac{|\nabla u|^2}2
   + W(u)
   - gu
   \Big).
\end{equation}
%
%
%
This leads to defining the {\sl energy space},
that will act as a phase space for our system:
\begin{equation}\label{defiX}
  \calX:=\big\{u\in V:~u\ge 0~\text{a.e.},~u^{1-\kappa} \in L^1(\Omega)\big\}.
\end{equation}
%
The space $\calX$
is endowed with the natural (graph) metric
\begin{equation}\label{defidX}
  \deriv_{\calX}(u_1,u_2):=
   \|u_1-u_2\|_V
   + \big\| u_1^{1-\kappa} - u_2^{1-\kappa} \big\|_{L^1(\Omega)},
\end{equation}
which is readily proved to be complete.
Given $\mu\in(0,\infty)$, we also define
\begin{equation}\label{defiXmu}
  \calX_\mu
   :=\big\{u\in \calX:~u\OO=\mu\big\},
\end{equation}
$(\cdot)\OO$ denoting here and below
the spatial average over $\Omega$.
Actually, integrating \eqref{tf1} in space
one readily sees that the
quantity $u\OO(t)$ is conserved
in time for any solution $u$.

Thanks to \eqref{hpgamma}, \eqref{hpg} and to
the above conservation property,
a direct computation
shows that, for some $\alpha_\mu,c_\mu,C_\mu>0$
(also depending on $g$),
there holds
\begin{equation}\label{Ecoerc}
  \alpha_\mu\Big(\|u\|^2_V
   + \big\|u^{1-\kappa}\big \|_{L^1(\Omega)} \Big)
   - C_\mu
  \le \calE(u) \le
   c_\mu \Big( 1 + \| u \|^2_V
   + \big\|u^{1-\kappa}\big \|_{L^1(\Omega)} \Big).
\end{equation}
for all $u\in \calX_\mu$. In particular, for a function $u$
of assigned spatial mean $\mu$, the finiteness
of the energy $\calE(u)$ corresponds exactly to the condition
$u\in \calX_\mu$.

\medskip

The above notation is sufficient to define the class
of weak solutions.
\bede\label{def:weaksol}
 A\/ {\rm weak solution} to problem~\eqref{tf1}-\eqref{neum-intro}
 is a couple $(u,w)$, with
 \begin{align}\label{regou}
   & u\in \LIV\cap\LDHD, \qquad
    \delta^{1/2} u \in \HUH,\\
   \label{regoenergy}
   & F(u) \in L^\infty(0,T;L^1(\Omega)),\\
  \label{regow}
   & w\in L^{5/4}(0,T;W^{1,{5/4}}(\Omega)),\\
  \label{regosergey}
   & b^{1/2}(u)\nabla w\in L^2(0,T;H),\\
  \label{regoult}
   & b(u)\nabla w\in L^{\frac{20}{10+s}}(Q), \qquad
    u_t\in L^{\frac{20}{10+s}}(0,T;(W^{1,\frac{20}{10-s}})^*(\Omega)),\\
  \label{regofu}
   & f(u)\in L^2(0,T;L^1(\Omega)) \cap L^{5/3}(Q).
 \end{align}
 such that the following relations hold
 a.e.~in~$(0,T)$:
 \begin{align}\label{tf1.2}
   & \duav{ u_t,\phi } + \duav{ b(u)\nabla w,\nabla \phi }=0
    \qquad\perogni\phi\in W^{1,\frac{20}{10-s}}(\Omega),\\
  \label{tf2.2}
   & \duav{ w, \psi} = \delta ( u_t,\psi )
    + \duav{ \nabla u,\nabla \psi} + \duav{ f(u) + \gamma(u) - g, \psi }
    \qquad\perogni\psi\in V, 
 \end{align}
 where $\langle \cdot,\cdot \rangle$ denote suitable
 duality pairings, and such that, in addition,
 \begin{equation}\label{init}
   u|_{t=0}=u_0,
   \quext{a.e.~in }\,\Omega.
 \end{equation}
\edde
\beos\label{enunc3d}
 The above regularity framework
 actually refers to the three-dimensional case. 
 If $d=2$, \eqref{regow}-\eqref{regofu} could be
 improved a bit (we omit the details).
\eddos
\noindent%
We can now state our basic existence result,
which can be considered as
a variant of the theorem proved in
\cite[Sec.~8]{Gr} (our assumptions are,
indeed, slightly different).
\bete\label{teoesi}
 Let us assume\/ \eqref{hpb}-\eqref{hpdelta}
 and let, for some $\mu\in(0,\infty)$,
 \begin{equation}\label{hpu0}
   u_0 \in \calX_\mu.
 \end{equation}
 Then, problem~\eqref{tf1}-\eqref{neum-intro}
 admits\/ {\rm at least} one weak solution. 
\ente
%


\section{Proof of Theorem~\ref{teoesi}}
\label{proof-esi}

We will give the proof in the case $d=3$ and just point out
some minor differences occurring in the two-dimensional case.


\subsection{Regularized problem}

First of all, we introduce a suitably approximated 
statement. Namely, given $\epsi\in(0,1)$, we set
\begin{align}\label{defibee}
  & b\ee(r) := b ( (r^2 + \epsi^a)^{1/2} ),\\
 \label{defifee}
  & f\ee(r) :=
   \begin{cases}
    f(r) & \quext{if }\,r\ge \epsi,\\
    \displaystyle f(\epsi)+f'(\epsi)(r-\epsi)
     =-\frac1{\epsi^\kappa}+\frac{\kappa}{\epsi^{\kappa+1}}(r-\epsi)
      & \quext{if }\,r< \epsi,
   \end{cases}
\end{align}
where $a>0$ will be chosen later on, and
it is intended that $f\ee$ is defined for all $r\in\RR$.
It is then worth noting that
\begin{equation}\label{feeprimo}
  f\ee'(r) :=
   \begin{cases}
    \displaystyle \frac\kappa{r^{\kappa+1}}
       & \quext{if }\,r\ge \epsi,\\[2mm]
     \displaystyle \frac\kappa{\epsi^{\kappa+1}}
       & \quext{if }\,r< \epsi.
   \end{cases}
\end{equation}
Moreover, setting
\begin{equation}\label{defiFee}
  F\ee(r) :=
   \frac1{\kappa-1}
    +\int_1^r f\ee(\tau)\,\ditau,
\end{equation}
we obtain that
\begin{equation}\label{Fee}
  F\ee(r) :=
   \begin{cases}
    \displaystyle \frac1{(\kappa-1)r^{\kappa-1}}
      & \quext{if }\,r\ge \epsi,\\[2mm]
    \displaystyle \frac1{(\kappa-1)\epsi^{\kappa-1}}
     - \frac1{\epsi^\kappa}(r-\epsi)
      +\frac\kappa{2\epsi^{\kappa+1}}(r-\epsi)^2
          & \quext{if }\,r< \epsi.
   \end{cases}
\end{equation}
At this point, we can consider the approximate statement
\begin{align}\label{tf1ee}
  & u\eet-\dive (b\ee(u\ee)\nabla w\ee)=0,\\
 \label{tf2ee}
  & w\ee = \delta u\eet - \Delta u\ee + f\ee(u\ee) + \gamma(u\ee) - g,
\end{align}
coupled with the initial conditions
and the no-flux boundary conditions.
Then, in analogy with \cite[Thms.~2 and 4]{EG}
(see also \cite[Thm.~2.1]{BB}), we have the following
existence result for approximate solutions.
\bete\label{teoepsi}
 Let us assume\/ \eqref{hpb}-\eqref{hpdelta}
 and let $b\ee$, $f\ee$ be specified by
 \eqref{defibee}-\eqref{defifee}. Let, in addition,
 $u\zee\in H^3(\Omega)$ be defined
 as the (unique) solution to the elliptic problem
 \begin{equation}\label{hpuzee}
   u\zee-\epsi^2\Delta u\zee=u_0,
    \qquad \dn u\zee|_{\partial \Omega}=0.
 \end{equation}
 Then, there exists at least one couple $(u\ee,w\ee)$
 with
 \begin{align}\label{regouee}
   & u\ee\in \LIV\cap L^2(0,T;H^3(\Omega)), \qquad
    \delta^{1/2} u\ee\in \HUH, \qquad \delta u\ee\in \HUV,    \\
  \label{regowee}
   & w\ee\in L^2(0,T;V),
 \end{align}
 satisfying\/ \eqref{tf1ee}-\eqref{tf2ee}
 a.e.~in~$\Omega\times(0,T)$, together with
 \begin{equation}\label{initee}
   u\ee|_{t=0}=u\zee,
    \quext{a.e.~in }\,\Omega.
 \end{equation}
 Moreover, in the case $\delta>0$, the couple $(u\ee,w\ee)$
 is unique.
\ente
\noindent%
{\bf Sketch of proof of Theorem~\ref{teoepsi}.}~~%
We can proceed by following very closely the proofs of
the quoted results \cite[Thms.~2 and 4]{EG},
\cite[Thm.~2.1]{BB}. For this reason, we will just give
some very brief highlights. Actually, the main difference here
is due to the the growth of $b\ee$ at infinity.
Nevertheless, one can of course truncate $b\ee$
near $\infty$ and replace it by some approximation
$b\een$ of at most linear growth, such that
$b\een$ suitably tends to $b\ee$ as $\nu\searrow 0$.
Then, existence for $\nu>0$ is proved similarly
as in \cite{EG,BB} and it remains to prove
suitable a-priori estimates uniform w.r.t.~$\nu$.
However, for the sake of simplicity we will omit
the $\nu$-approximation and rather perform
formal estimates on the $\epsi$-solution in order
to show that it fulfills regularity properties
\eqref{regouee}-\eqref{regowee}.

Thus, we can firstly perform the {\sl energy}\/
and {\sl entropy}\/ estimate, as below.
Using that $b\ee\ge \epsi^{a s/2}$
and the Lipschitz continuity of $f\ee$,
we then obtain
\begin{align}\label{stuee-1}
  & \|u\ee\|_{\LIV\cap\LDHD}
   + \delta^{1/2} \|u\ee\|_{\HUH}\le c\ee,\\
 \label{stwee-1}
  & \|b\ee^{1/2}\nabla w\ee\|_{\LDH}
   +\|w\ee\|_{\LDV}
  \le c\ee.
\end{align}
Next, we test \eqref{tf2ee} by $\Delta^2 u\ee$.
Using \eqref{hpuzee}, \eqref{stuee-1}-\eqref{stwee-1},
\eqref{hpg} and the global Lipschitz continuity
of $f\ee$ and $\gamma$, it is then not difficult to obtain
\begin{equation}\label{stuee-2}
  \|u\ee\|_{L^2(0,T;H^3(\Omega))}
   + \delta^{1/2} \|u\ee\|_{L^\infty(0,T;H^2(\Omega))}\le c\ee.
\end{equation}
In case $\delta>0$, we can also test \eqref{tf2ee} by
$-\Delta u\eet$, that yields
\begin{equation}\label{stuee-3}
  \delta \|u\ee\|_{H^1(0,T;V)}\le c\ee.
\end{equation}
This gives all desired regularity properties.
Finally, to prove uniqueness in the case $\delta>0$, we
can test the difference of \eqref{tf2ee} by
the difference of the $u\eet$ and
the difference of \eqref{tf1ee} by the difference
of $w\eet$. The details are left to the 
reader.~\dimbox


\subsection{A priori estimates}

We now aim to obtain a number of a priori bounds,
uniform in $\epsi$, with the purpose of removing the
approximation.

\smallskip

\noindent%
{\bf Energy estimate.}~~%
We test \eqref{tf1ee} by $w\ee$, \eqref{tf2ee} by
$u\eet$, and sum the results. We then obtain
\begin{equation}\label{conto11}
  \ddt\calE\ee(u\ee)
   +\delta\|u\eet\|^2
   +\io b\ee(u\ee)|\nabla w\ee|^2
  =0,
\end{equation}
where
\begin{equation}\label{defiEee}
  \calE\ee(u):=\io\Big(
   \frac{|\nabla u|^2}2
   +F\ee(u)+\Gamma(u)-gu
    \Big).
\end{equation}
\beos\label{pocarego}
 We point out that, even at the approximate level,
 this formal estimate is not completely justified
 in the case $\delta=0$. Actually, $w\ee$ is
 (only) in $\LDV$, while \eqref{tf1ee} is not
 an equation in $\LDVp$ since $b\ee$ grows
 fast at infinity. However it is clear that, performing
 a truncation of $b\ee$ and then passing to
 the limit, the estimate could be justified.
\eddos
\noindent%
Then, we integrate \eqref{conto11} in time and
notice that, by \eqref{Fee},
\begin{equation}\label{Fee0}
  \io F\ee(u\zee) 
   \le \io F(u\zee)
   \le \io F(u_0) + \big( f(u\zee), u\zee - u_0 \big)
   \le \io F(u_0),
\end{equation}
the latter inequality following from \eqref{hpuzee} and 
Green's formula. Thus, owing to \eqref{hpu0},
\eqref{conto11} gives
\begin{align}\label{st11}
  & \|u\ee\|_{\LIV}
   +\delta^{1/2}\|u\eet\|_{\LDH}\le c,\\
 \label{st12}
  & \|F\ee(u\ee)\|_{L^\infty(0,T;L^1(\Omega))}
   \le c,\\
 \label{st13}
  & \|b\ee^{1/2}(u\ee)\nabla w\ee\|_{L^2(0,T;H)}
   \le c.
\end{align}
Here and below, $c$ denotes a positive constant,
independent of $\epsi$ and of time, whose value
may vary even inside a single row. We will use
the letter $\alpha$ to denote constants used in
estimates from below.

\smallskip

\noindent%
{\bf Entropy estimate.}~~%
Let us define, for $r\in(0,\infty)$,
the {\sl entropy}~$M$, by setting
\begin{equation}\label{defim}
  m(r):=\int_1^r \frac{\ditau}{b(\tau)}, \qquad
   M(r):=\int_1^r m(\tau)\,\ditau.
\end{equation}
Clearly, $M$ is a convex function
that grows at most like $r$ for $r\sim\infty$.
Moreover, let us introduce its approximate version
by taking, for $r\in\RR$,
\begin{equation}\label{defimee}
  m\ee(r):=\int_1^r \frac{\ditau}{b\ee(\tau)}, \qquad
   M\ee(r):=\int_1^r m\ee(\tau)\,\ditau.
\end{equation}
Clearly, $M\ee$ is a convex function
such that $M\ee\le M$ a.e.~in $(0,\infty)$.
Moreover, $m\ee$, $M\ee$ tend
to $m$, $M$, respectively, uniformly on compact
sets of $(0,\infty)$.

Then, we can test \eqref{tf1ee} by
$m\ee(u\ee)$, \eqref{tf2ee} by $-\Delta u\ee$,
and sum the results. We deduce
\begin{equation}\label{st21}
  \ddt \Big(\io M\ee(u\ee)
   +\frac\delta2\|\nabla u\ee\|^2\Big)
   +\|\Delta u\ee\|^2
   +\io \big(f\ee'(u\ee)+\gamma'(u\ee)\big)|\nabla u\ee|^2
   +\big(g,\Delta u\ee\big)
  =0,
\end{equation}
and we have to control some terms. Firstly,
by \eqref{hpgamma} and H\"older's inequality,
we have
\begin{equation}\label{st22}
  \bigg|\io \gamma'(u\ee)|\nabla u\ee|^2
   +\big(g,\Delta u\ee\big) \bigg|
   \le c\big(1+\|\nabla u\ee\|^2\big)
    +\frac12\|\Delta u\ee\|^2.
\end{equation}
We now observe that, for all $\epsi\in(0,1)$,
it is $M\ee(u\zee) \le M(u_0)$ (to prove this,
proceed as in \eqref{Fee0}). Moreover, we have
that $M(u_0)\in L^1(\Omega)$. Actually,
\eqref{hpu0} entails $u_0^{1-\kappa}\in L^1(\Omega)$
and $M(r)$ grows no faster than
$r^{2-s}$ in the neighbourhood
of $0$, which is good since we assumed
$\kappa\ge s+1$ (cf.~\eqref{hpf}).
Thus, integrating \eqref{st22}
in time we arrive at
\begin{align}\label{st23}
  & \|u\ee\|_{\LDHD} \le c,\\
 \label{st24}
  & \iTo f\ee'(u\ee)|\nabla u\ee|^2\le c.
\end{align}

\smallskip

\noindent%
{\bf Control of the nonlinear terms.}~~%
Here, our aim is to derive $\epsi$-uniform bounds
in order to pass to the limit in
the terms $b\ee(u\ee)$ and $f\ee(u\ee)$.

First of all, by \eqref{st11}, \eqref{st23} and
interpolation, we obtain, if $d=3$,
\begin{equation}\label{st31}
  \|u\ee\|_{L^{10}(\Omega\times(0,T))} \le c,
\end{equation}
whereas for $d=2$ we have instead
\begin{equation}\label{st31-2d}
  \|u\ee\|_{L^{\ppp}(\Omega\times(0,T))} \le c_\ppp
   \quad\perogni \ppp\in[1,\infty).
\end{equation}
Let us now set
\begin{align}\label{defQee-2}
  & Q\ee^1:=\big\{(x,t)\in Q:~u\ee(x,t)<\epsi\big\},\\[1mm]
 \label{defQee-1}
  & Q\ee^2:=\big\{(x,t)\in Q:~\epsi \le u\ee(x,t)\le 1\big\},\\[1mm]
 \label{defQee}
  & Q\ee^3:=\big\{(x,t)\in Q:~u\ee(x,t)>1\big\}
\end{align}
and notice that, by \eqref{Fee} and \eqref{st12},
\begin{equation}\label{measQ1}
  |\Omega\ee^1(t)|\le c\epsi^{\kappa-1}
   \quext{for a.e.~}\,t\in(0,T).
\end{equation}
Here and below, $\Omega\ee^i(t)$, $i\in\{1,2,3\}$,
is the section of $Q\ee^i$ at the
generic time $t\in(0,T)$.
Also, it is obvious that
\begin{equation}\label{coQ3}
  \| f\ee(u\ee) \|_{L^\infty(Q\ee^3)}\le c.
\end{equation}
%
%
%
%
%
%
%
Next, we remark (cf.~\eqref{defibee}) that
\begin{equation}\label{buinv}
  \frac1{b\ee(r)}\le
   \begin{cases}
    \displaystyle \frac{c}{r^s+\beta r^n}
       & \quext{if }\,r\ge \epsi,\\[2mm]
     \displaystyle \frac{c}{\epsi^{a s/2}}
       & \quext{if }\,r< \epsi.
   \end{cases}
\end{equation}
Hence, using \eqref{Fee}, \eqref{st12},
\eqref{measQ1} and the condition $\kappa\ge s+1$ 
(cf.~\eqref{hpf}), taking $a$ small enough in 
the definition
\eqref{defibee} of $b\ee$ we obtain
\begin{equation}\label{stbee1}
  \big\| b\ee^{-1}(u\ee) \big\|_{L^\infty(0,T;L^1(\Omega))}
   \le c.
\end{equation}
Setting now
\begin{equation}\label{defivee}
  z\ee(x,t):=\max\big\{u\ee(x,t),\epsi\big\}, \qquad
   v\ee(x,t):=z\ee(x,t)^{\frac{1-\kappa}2},
\end{equation}
for a.e.~$t\in(0,t)$ we have
\begin{align}\no
  \| v\ee(t) \|_{L^6(\Omega)}^2
   & \le c \big( \| v\ee(t) \|^2
    + \| \nabla v\ee(t) \|^2 \big)\\
  \label{contovee}
   & \le c \int_{\Omega\ee^1(t)} \frac1{\epsi^{\kappa-1}}
    + c \int_{\Omega\ee^2(t)\cup \Omega\ee^3(t)} \frac1{u\ee^{\kappa-1}}
    + c \int_{\Omega\ee^2(t)\cup \Omega\ee^3(t)}
     \frac{|\nabla u\ee|^2}{u\ee^{\kappa+1}}
\end{align}
(here and below, we compute the exponents referring
to the case $d=3$; for $d=2$, they can be improved,
of course).
Hence, using \eqref{st12}, \eqref{st21}
and \eqref{measQ1}, and
integrating in time, it is not difficult
to arrive at
\begin{equation}\label{stvee1}
  \| v\ee \|_{L^2(0,T;L^6(\Omega))}
   \le c,
\end{equation}
whence, recalling \eqref{buinv} and possibly choosing
a smaller $a$,
\begin{equation}\label{stbee2}
  \big\| b\ee^{-1}(u\ee) \big\|_{L^1(0,T;L^3(\Omega))}
   \le c.
\end{equation}
Then, using \eqref{st13} and either \eqref{stbee1} or
\eqref{stbee2}, we also have
\begin{equation}\label{stwee1}
  \| \nabla w\ee \|_{L^2(0,T;L^1(\Omega))}
   + \| \nabla w\ee \|_{L^1(0,T;L^{3/2}(\Omega))}
   \le c,
\end{equation}
whence, by interpolation,
\begin{equation}\label{stwee2}
  \| \nabla w\ee \|_{L^{5/4}(Q)}
   + \| \nabla w\ee \|_{L^{5/3}(0,T;L^{15/14}(\Omega))}
   \le c.
\end{equation}
We are now ready to
give an estimate of the term $f\ee(u\ee)$.
To do this, we first test \eqref{tf2ee} by
$u\ee-\mu$, to obtain
\begin{equation}\label{contomu1}
  \mu \io | f\ee(u\ee) |
   =  \io \big( - \delta u\eet + \Delta u\ee - \gamma(u\ee) + g \big)
     (u\ee - \mu)
   + \io w\ee ( u\ee - \mu )
   - \io f\ee(u\ee) u\ee
\end{equation}
and the terms on the \rhs\ are treated as follows:
\begin{align}\no
  \io \big( - \delta u\eet + \Delta u\ee - \gamma(u\ee) + g \big)
     (\mu - u\ee)
   & \le c \big\| - \delta u\eet + \Delta u\ee - \gamma(u\ee) + g \big\|
    \| \mu - u\ee \| \\
 \label{contomu2}
  & \le c \big\| - \delta u\eet + \Delta u\ee - \gamma(u\ee) + g \big\|
    =: \eta_1,
\end{align}
thanks to \eqref{st11}, where $\|\eta_1\|_{L^2(0,T)}\le c$ by
\eqref{hpgamma}, \eqref{hpg} and \eqref{st23}. Next,
\begin{equation}\label{contomu3}
  \io w\ee ( u\ee - \mu )
  = \io \big( w\ee - (w\ee)\OO \big) ( \mu - u\ee )
   \le c \| \nabla w\ee \| 
   =: \eta_2,
\end{equation}
with $\|\eta_2\|_{L^2(0,T)}\le c$ by
\eqref{stwee1}. Finally, due
to \eqref{defifee} it is clear that
\begin{equation}\label{contomu4}
  - \io f\ee(u\ee) u\ee
   \le \frac\mu2 \io | f\ee(u\ee) |
   + c_\mu. 
\end{equation}
Squaring \eqref{contomu1}, integrating in time,
and using \eqref{contomu2}-\eqref{contomu4}, 
we easily arrive at
\begin{equation}\label{stmu}
  \| f\ee(u\ee) \|_{L^2(0,T;L^1(\Omega))}
   \le c. 
\end{equation}

\smallskip

Setting now $\phi\ee:=-|f\ee(r)|^{2/3}$, we
test \eqref{tf2ee} by $\phi\ee(u\ee)-(\phi\ee(u\ee))\OO$,
and integrate in space and time. This gives
\begin{align}\no
  & \iTo f\ee(u\ee) \phi\ee(u\ee)
   + \iTT \big( \delta u\eet + \gamma(u\ee) - g ,
      \phi\ee(u\ee)-(\phi\ee(u\ee))\OO \big)\\
  \label{conto51}
  & \mbox{}~~~~~
   = \iTT \big( w\ee - (w\ee)\OO, \phi\ee(u\ee) \big)
    + \iTT \Big( (\phi\ee(u\ee))\OO \io f\ee(u\ee) \Big).
\end{align}
Let us first notice that the first term on the \lhs\ gives
\begin{equation}\label{conto52}
  \iTo f\ee(u\ee) \phi\ee(u\ee)
   = \|f\ee(u\ee)\|^{5/3}_{L^{5/3}(Q)}.
\end{equation}
Next, by H\"older's inequality, we obtain
\begin{align}\no
  & \bigg| \iTT \big( \delta u\eet + \gamma(u\ee) - g ,
      \phi\ee(u\ee)-(\phi\ee(u\ee))\OO \big)\bigg|\\
 \label{conto54}
   & \mbox{}~~~~~
  \le \|\delta u\eet + \gamma(u\ee) + g\|_{L^2(Q)}^2
   + c \| \phi\ee(u\ee) \|_{L^2(Q)}^2
  \le \sigma \|f\ee(u\ee)\|_{L^{5/3}(Q)}^{5/3}
   +c_\sigma.
\end{align}
where $\sigma$ is a small constant to be chosen below
and $c_\sigma>0$ depends on $\sigma$. Notice that
\eqref{st11}, \eqref{hpgamma} and \eqref{hpg}
have been used here.

It then remains to control the terms on the
\rhs\ of \eqref{conto51}. As far as the first
one is concerned, recalling \eqref{coQ3},
we notice that, a.e.~in $(0,T)$,
\begin{align}\no
  \big( w\ee - (w\ee)\OO, \phi\ee(u\ee) \big)
   & \le \|\phi\ee(u\ee)\|_{L^{5/2}(Q)}
    \|w\ee - (w\ee)\OO\|_{L^{5/3}(Q)}\\
 \label{conto55}
  & \le \sigma \|f\ee(u\ee)\|_{L^{5/3}(\Omega)}^{5/3}
   + c_\sigma \|\nabla w\ee\|_{L^{15/14}(\Omega)}^{5/3}.
\end{align}
Finally, let us estimate the latter term in \eqref{conto51}.
Actually, it is clear that
\begin{equation}\label{conto57b}
  \iTT (\phi\ee(u\ee))\OO \io f\ee(u\ee)
   \le c  \|f\ee(u\ee)\|_{L^2(0,T;L^1(\Omega))}^2
    + c 
   \le c,
\end{equation}
the latter inequality following from \eqref{stmu}.
%
%
%
%
%
Collecting now \eqref{conto51}-\eqref{conto57b}
and recalling \eqref{coQ3}, we finally arrive at
\begin{equation}\label{st51}
  \| f\ee(u\ee) \|_{L^{5/3}(Q)}\le c.
\end{equation}
%


\subsection{Passage to the limit}
\label{pass-lim}

We will just consider, for brevity, the 
case $d=3$.
For simplicity of notation, let us set
$\zeta\ee:=b\ee^{1/2}(u\ee)\nabla w\ee$
and let $\| \cdot \|_p$ denote
the norm in the space $L^p(Q)$.
Then, by \eqref{st13},
\begin{equation}\label{st71}
  \| \zeta\ee \|_{2}\le c.
\end{equation}
Moreover, being $s<10$, \eqref{st31}
guarantees that
\begin{equation}\label{st51.3}
  \big\| b\ee(u\ee) \nabla w\ee \big\|_q
   = \big\| b\ee^{1/2}(u\ee) \zeta\ee \big\|_q
   \le c \quext{for some }\,q>1.
\end{equation}
Hence, by comparison in \eqref{tf1ee},
\begin{equation}\label{st73}
  \| u\eet \|_{L^q(0,T;W^{-1,q}(\Omega))}\le c
   \quext{for some }\,q>1
\end{equation}
(of course, if $\delta>0$ we have much more, but
we want to deal with the most general case here).
Consequently, using \eqref{st11}, \eqref{st31},
the Aubin-Lions lemma, and Lebesgue's theorem,
\begin{equation}\label{co81}
  u\ee\to u \quext{strongly in }\,
   L^q(\Omega\times(0,T)) \quad \perogni q\in[1,10)
\end{equation}
(here and below, all convergence relations are to be intended
up to the extraction of non-relabelled subsequences).
Let us now notice that
\begin{equation}\label{cocomp}
  b\ee\to \bar b
   \quext{uniformly on compact subsets of }\,\RR,
\end{equation}
where $\bar b$ denotes the even extension of $b$ to
$\RR$. Then, by \eqref{co81}, $s<10$ in \eqref{hpb},
and Lebesgue's theorem again, we obtain
\begin{equation}\label{co82}
  b\ee(u\ee)\to \bar b(u) \quext{strongly in }\,
   L^q(\Omega\times(0,T)) \quext{for some }\,q>1.
\end{equation}
Analogously, by \eqref{st51},
\begin{equation}\label{co83}
  f\ee(u\ee)\to f(u) \quext{strongly in }\,
   L^q(\Omega\times(0,T)) \quext{for all }\,q\in[1,5/3),
\end{equation}
whence, in particular, the limit $u$ is a.e.~nonnegative
and we can replace $\bar b$ with $b$ in \eqref{co82}.

Our next aim is to pass to the limit in the product
$b\ee(u\ee)\nabla w\ee$. To do this, we first notice that,
by \eqref{co81} and $s<10$,
\begin{equation}\label{st51.3a}
  b\ee(u\ee)^{1/5}
   \to b(u)^{1/5} \quext{strongly in }\,L^q(Q)
   \quext{for some }\,q>5.
\end{equation}
Thus, using also \eqref{stwee2} we arrive at
\begin{equation}\label{st51.3b}
  b\ee(u\ee)^{1/5}\nabla w\ee
   \to b(u)^{1/5}\nabla w
   \quext{weakly in }\,L^q(Q)
   \quext{for some }\,q>1.
\end{equation}
Next, interpolating between \eqref{stbee1}
and \eqref{stbee2}, we get
\begin{equation}\label{st51.2}
  \big\| b\ee(u\ee)^{-1} \big\|_{5/3}\le c.
\end{equation}
Using \eqref{st51.2} and \eqref{st71},
we then obtain
\begin{equation}\label{st51.4}
  \big\| b\ee(u\ee)^{1/5}\nabla w\ee \big\|_{25/17}
   \le \| \zeta\ee \|_2 \big\| b\ee(u\ee)^{-3/10} \big\|_{50/9}
   \le c.
\end{equation}
Hence, by \eqref{st51.4}, \eqref{st51.3a}
and \eqref{st51.3b},
\begin{equation}\label{st51.4b}
  b\ee(u\ee)^{2/5}\nabla w\ee
   = b\ee(u\ee)^{1/5} \big( b\ee(u\ee)^{1/5}\nabla w\ee \big)
   \to b(u)^{2/5}\nabla w
  \quext{weakly in }\,L^{25/22}(Q).
\end{equation}
Then, writing $b\ee(u\ee)^{2/5}\nabla w\ee$
as $\zeta\ee b\ee(u\ee)^{-1/10}$ and using once more
\eqref{st71} and \eqref{st51.2}, it is clear that
the exponent $25/22$ in \eqref{st51.4b} can be improved.
Thus, iterating the above procedure and using
once more \eqref{st71}, it is not difficult to
arrive at
\begin{equation}\label{st51.6}
  \zeta\ee \to b(u)^{1/2}\nabla w
   \quext{weakly in }\,L^2(Q),
\end{equation}
whence, using once more \eqref{co82}
(and computing explicitly the exponent),
we finally obtain
\begin{equation}\label{st51.7}
  b\ee(u\ee)\nabla w\ee \to b(u)\nabla w
   \quext{weakly in }\,L^{\frac{20}{10+s}}(Q).
\end{equation}
Thus, we can take the limit of all terms in
\eqref{tf1ee}-\eqref{tf2ee} and get back
\eqref{tf1.2}-\eqref{tf2.2}, where the
properties required to test functions
depend of course on the regularity
conditions proved above. Indeed,
\eqref{regou}-\eqref{regofu}
follow as a direct byproduct of the 
procedure. The proof of
Theorem~\ref{teoesi} is concluded.


\section{Weak trajectory attractors}
\label{sec:weakattr}

In this section, we construct the so-called weak trajectory attractor 
for problem \eqref{tf1}-\eqref{tf2}.
To avoid technicalities, we will limit ourselves
to deal with the (more degenerate) case $\beta=0$, i.e.,
$b(u)=u^s$, for $s\in [0,10)$ (cf.~\eqref{hpb}).

We start by proving a dissipativity result holding for the
weak solutions constructed in the proof of Theorem~\ref{teoesi}.
\bete\label{teodiss}
 Let the assumptions of\/ {\rm Theorem~\ref{teoesi}}
 hold, with $\beta=0$, and let
 \begin{equation}\label{defiE0}
   \EE_0 := \calE(u_0),
 \end{equation}
 that is finite thanks to \eqref{hpu0}. Then,
 there exist a weak solution $(u,w)$ and a monotone function
 $Q:[0,\infty)\to [0,\infty)$ such that
 \begin{equation}\label{energyest}
   \calE(u(t))
     + \itt \big( \delta \| u_t \|^2
      + \| b^{1/2}(u) \nabla w \|^2 \big)
    \le Q(\EE_0)
    \quad\perogni t\ge 0.
 \end{equation}
 More precisely, there exists
 a set $\calB_0$ bounded with respect to
 the metric \eqref{defidX}, such that,
 for any $\deriv_{\calX}$-bounded set
 $B\subset\calX_\mu$, there exist a time $T_B\ge 0$ 
 such that for any initial datum $u_0\in B$ there
 exists at least one weak solution $u$ starting
 from $u_0$ and such that $u(t)\in\calB_0$ for all 
 $t\ge T_B$.
\ente
\beos
 Notice that the above result does not claim that
 dissipativity holds in the whole class of weak solutions,
 but just that from any admissible initial datum there 
 starts (at least) one weak solution in the dissipative
 class (cf.~Remark~\ref{RemW.bounded} for further considerations).
\eddos
\begin{proof}
We integrate \eqref{conto11} between
$0$ and an arbitrary $t>0$. This gives the $\epsi$-equivalent
of \eqref{energyest}. Then, we take the $\liminf$ with respect
to $\epsi\searrow 0$ and use estimates \eqref{st11},
\eqref{co81}, \eqref{st51.6}, the fact that $F\ee$ converges
to $F$ uniformly on compact sets of $(0,\infty)$, Fatou's
Lemma, and the lower semicontinuity of norms with
respect to weak or weak star convergences (of course,
we will get an energy {\sl inequality}, and not necessarily
an equality, in this way). Relation~\eqref{energyest} is proved.

To show dissipativity, we start considering
the case $\kappa>s+1$. Then, it is convenient
to rewrite the energy inequality in the differential form:
\begin{equation}\label{energylim}
  \ddt\calE(u)
   +\delta \| u_t \|^2
   +\io b(u)|\nabla w|^2
  \le 0,
\end{equation}
for a.e.~$t>0$. We then notice that, 
setting $z:=u^{-1}$, we have, from \eqref{Ecoerc},
\begin{equation}\label{propE}
  \sigma \big( \| z^{\kappa-1} \|_{L^1(\Omega)}
   + \| \nabla u \|^2 \big) - c
  \le \calE(u)
   \le \sigma^{-1} \big( \mu \| z \|_{L^{\kappa}(\Omega)}^{\kappa}
   + \| \nabla u \|^2 + 1 \big),
\end{equation}
where $\sigma\in(0,1)$ depends in particular on $\mu=u\OO$.

Now, similarly with \eqref{contomu1}, we
test \eqref{tf2} by $u-\mu$, obtaining
\begin{align}\no
  \| \nabla u \|^2
   + \mu \| z \|_{L^\kappa(\Omega)}^\kappa
  & \le \io (- \delta u_t - \gamma(u) + g ) (u - \mu)
   + \| z^{\kappa-1} \|_{L^1(\Omega)}
   + \io w (u - \mu)\\
 \no
  & \le \frac12 \| \nabla u \|^2
   + c \big (1 + \delta^2 \| u_t \|^2 \big)
   + \frac\mu2\| z \|_{L^{\kappa}(\Omega)}^{\kappa}
   + c_\mu
   + c \| w - w\OO \|_{L^{3/2}(\Omega)}^2\\
 \label{contol2l1-1prev}
  & \le \frac12 \| \nabla u \|^2
   + c \big (1 + \delta^2 \| u_t \|^2 \big)
   + \frac\mu2\| z \|_{L^{\kappa}(\Omega)}^{\kappa}
   + c_\mu
   + c \| \nabla w \|_{L^1(\Omega)}^2.
\end{align}
The last term can be controlled this way:
\begin{equation}\label{contoen11}
  \| \nabla w \|_{L^1(\Omega)}^2
   \le c \| z^s \|_{L^1(\Omega)} \io b(u) | \nabla w |^2.
\end{equation}
Next, being $\kappa > s+1$, using the first inequality in
\eqref{propE} and Jensen's inequality, we obtain,
for suitable positive constants $\alpha$,
\begin{equation}\label{contoen12}
  \big(\calE(u)\big)^{\frac{s}{\kappa-1}}
   \ge \alpha \| z^s \|_{L^1(\Omega)}
   \ge \alpha (z\OO)^s \ge \alpha (u\OO)^{-s}
   = \alpha \mu^{-s}.
\end{equation}
From \eqref{propE}-\eqref{contoen12}, we then obtain,
for a suitable $c_*>0$,
\begin{equation}\label{contoen13}
  \calE(u)
   \le c_* \Big (1 + \delta^2 \| u_t \|^2
   + \big(\calE(u)\big)^{\frac{s}{\kappa-1}}
     \io b(u)|\nabla w|^2 \Big).
\end{equation}
Thus, assuming $\calE\ge1$, which is of course not restrictive,
we can divide by $\calE^{s/(\kappa-1)}$ to obtain
\begin{equation}\label{contoen13bis}
  \big(\calE(u)\big)^{\frac{\kappa-1-s}{\kappa-1}}
   \le c_* \Big (1 + \delta^2 \| u_t \|^2
   + \io b(u)|\nabla w|^2 \Big).
\end{equation}
Taking the $(1-\epsilon)$-power for $\epsilon\in(0,1)$,
summing to \eqref{energylim}, and applying
Young's inequality, we then have, for some
$\alpha>0$,
\begin{equation}\label{contoen16}
  \ddt\calE(u)
   + \alpha \Big( \big(\calE(u)\big)^{\frac{(\kappa-1-s)(1-\epsilon)}{\kappa-1}}
   + \delta \| u_t \|^2
   + \io b(u)|\nabla w|^2 \Big) \le c_\epsilon,
\end{equation}
whence the thesis follows by integrating in time and applying
the comparison principle for ODE's. More precisely, we also
have a quantitative decay estimate for the energy.

In the case $\kappa=s+1$, we can say a little bit less, but
dissipativity still holds. Actually,
we can repeat the above procedure up to \eqref{contoen13}.
Then, we notice that \eqref{energylim} implies
in particular that
\begin{equation}\label{dissinte}
   \int_0^\infty \Big(
    \delta \| u_t \|^2
   + \io b(u)|\nabla w|^2 \Big)
  \le \EE_0<\infty.
\end{equation}
Consequently, there exists at least one
time $T_*$ such that
\begin{equation}\label{Tstar}
  T_*\in [0,2c_*\EE_0], \quext{and}~~
   \delta \| u_t(T_*) \|^2
   + \io b(u(T_*))|\nabla w(T_*)|^2
  \le \frac1{2c_*}.
\end{equation}
Substituting in \eqref{contoen13}, we then have
\begin{equation}\label{contoen13diss}
  \calE(u(T_*)) \le 2 c_* + \delta
   =:C^*.
\end{equation}
Then, we obtain $\calE(u(t))\le C^*$ for any
$t\ge T_*$ simply observing that, by \eqref{energylim},
$\calE$ is nonincreasing.
\end{proof}
\noindent%
As a next step (see \cite{CV} for more details), 
we need to rewrite the dissipative estimate in such a way
that, on the one hand, we will be able to control all the norms 
which are necessary to pass to the weak limit
in the space of solutions of the problem considered 
(and verify that the limit function 
is again a solution) and, on the other hand, be sure that the corresponding 
trajectory phase space will be translation invariant.

The following lemma improves the dissipative estimate \eqref{energyest} 
by adding the terms controlled by the entropy estimate.
\bele\label{LemW.energy} 
 Let the assumptions of\/ {\rm Theorem \ref{teodiss}} hold. Then, 
 there exists a solution $(u,w)$ of problem\/ 
 \eqref{tf1}-\eqref{tf2} which satisfies the following estimate:
 \begin{align}\no
  & \mathcal E(u(T))
   + \int_T^{T+1}\delta\|u_t(t)\|^2
   + \|b^{1/2}(u(t))\Nx w(t)\|^2
   + \|\Dx u(t)\|^2
   + \big(f'(u(t))\Nx u(t),\Nx u(t)\big)\,\dit\\
  \label{W.en}
  & \mbox{}~~~~~~~~~~
   \le Q(\mathcal E(u(0)))e^{-\alpha T}+C_*,
   \quad \perogni T\ge 0.
 \end{align}
 where the positive constants $\alpha$ and $C_*$ and the monotone 
 function $Q$ are independent of $t$ and of 
 the concrete choice of the solution $u$.
\enle
\begin{proof}
As a consequence of Theorem~\ref{teodiss}, it is clear that \eqref{W.en}
holds, for suitable $Q$ and $\alpha$, without the last two terms
in the integral on the \lhs. To control these terms,
it is sufficient, in the case $\kappa>s+1$, to sum~\eqref{st21} to
\eqref{contoen16} in the preceding proof. In the case $\kappa=s+1$, 
we know that there exists $T_*=T_*(\EE_0)$ such that the energy is 
smaller than some constant $C^*$ independent of the initial data for
any $T\ge T_*$ (cf.~\eqref{contoen13diss}). Then, it is sufficient
to integrate \eqref{st21} over $(T,T+1)$ for $T\ge T_*$
to get
\begin{align}\no
  & \int_T^{T+1}
   \|\Delta u(t)\|^2
   + \big( f'(u(t)) \nabla u(t), \nabla u(t) \big)\,\dit\\
 \label{st21inte}
  & \mbox{}~~~~~
  \le c + \io M(u(T))
   + \frac\delta2 \|\nabla u(T)\|^2
  \le Q(\calE(T)) 
  \le Q(C^*),
\end{align}
as desired. Actually, it is clear that the energy
controls from above the terms $M(u)$ and $\|\Nx u\|$.
Being pedantic, all these estimates should be done on the 
level of approximations $u_\eb$ with passing to the
limit after that; we directly performed the estimate
on $u$ just for brevity.
\end{proof}
\beos\label{RemW.rest}
 Note also that all the norms involved into our definition of a weak solution 
 (see Theorem \ref{teoesi}) are under control if we assume that the weak 
 solution satisfies \eqref{W.en}. This fact, which can be verified exactly as in 
 the proof of Theorem \ref{teoesi}, is crucial in order to be able to 
 pass to the weak limit on the space of weak solutions and establish that the 
 absorbing set for the trajectory semigroup is indeed closed, see below.
\eddos
\noindent%
We are now able to define the trajectory phase space and trajectory dynamical system associated 
with problem \eqref{tf1}-\eqref{tf2}.
\bede
 Let $\mathcal K_+\subset L^\infty(\R_+,\mathcal E)$ be the set of all solutions $u$ 
 of problem \eqref{tf1}-\eqref{tf2} belonging to the class \eqref{regou}-\eqref{regofu} 
 which satisfy the following analogue of \eqref{W.en}:
 \begin{align}\no
  & \mathcal E(u(T))
   + \int_T^{T+1}\delta\|u_t(t)\|^2
   + \|b^{1/2}(u(t))\Nx w(t)\|^2
   + \|\Dx u(t)\|^2
   + (f'(u(t))\Nx u(t),\Nx u(t))\,\dit\\
  \label{W.en1}
   & \mbox{}~~~~~~~~~~
   \le C_ue^{-\alpha T}+C_*,
   \quad \perogni T\ge 0.
 \end{align}
 for some constant $C_u$ depending on the solution $u$. Then, 
 the shift semigroup $T(h)$, $h\ge0$, acts on $\mathcal K_+$:
 \begin{equation}\label{W.sem}
   T(h):\mathcal K_+\to\mathcal K_+,
    \qquad (T(h)u)(t):=u(t+h).
 \end{equation}
 We will refer below to $\mathcal K_+$ and 
 $T(h):\mathcal K_+\to\mathcal K_+$ as a trajectory phase space and 
 trajectory dynamical system associated with problem \eqref{tf1}-\eqref{tf2} respectively.
\edde
Furthermore, in order to be able to introduce the attractor of the 
trajectory dynamical system, we need to specify the topology on 
$\mathcal K_+$ as well as the class of bounded sets.
\bede\label{DefW.bound} 
 We endow the set $\mathcal K_+$ with the topology induced by the embedding
 $\mathcal K_+\subset\Theta_+^{weak}:=
 [L^\infty_{\loc}(\R_+,H^1(\Omega))\cap L^2_{\loc}(\R_+,H^1(\Omega))]^{w^*}$,
 where $w^*$ stands for the weak-star topology, and will refer to it as a 
 weak topology on the trajectory phase space $\mathcal K_+$. 
 \par
 A set $B\subset \mathcal K_+$ will be called bounded if inequality\/
 \eqref{W.en1} holds uniformly with respect to all $u\in B$, i.e., if
 \begin{equation}\label{W.bound}
   C_B:=\sup_{u\in B}C_u<\infty.
 \end{equation}
\edde
\beos\label{RemW.bounded} 
 As usual (see \cite{CV} for the details), under the general assumptions of 
 Theorem \ref{teodiss}, we know neither the fact that any weak solution 
 of problem \eqref{tf1}-\eqref{tf2} satisfies
 the energy inequality \eqref{W.en1} nor that the constant $C_u$ in \eqref{W.en1} 
 can be expressed in terms of $\mathcal E(u(0))$. Actually, it may be possible to construct a 
 solution $u$ which satisfies \eqref{W.en} for the initial moment $T=0$ only and be
 unable to verify its analogue for other initial times. By this reason, 
 attempting to replace \eqref{W.en1} by \eqref{W.en} in the definition of 
 the trajectory phase space $\mathcal K_+$, we lose the translation invariance 
 $T(h)\mathcal K_+\subset \mathcal K_+$ which is crucial for the attractors theory. 
 However, as we will see below, under the more restrictive assumptions of 
 Theorem~\ref{teointerm}, the answer to both the questions posed above 
 is positive. So, in that case, every reasonably defined weak solution
 satisfies \eqref{W.en} and the boundedness condition is equivalent
 to the boundedness of $u(0)$ in the energy space.
\eddos
\noindent%
Finally, we are now able to introduce the trajectory 
attractor for problem \eqref{tf1}-\eqref{tf2}.
\bede\label 
 A set $\mathcal A^{tr}\subset\mathcal K_+$ is a (weak) trajectory attractor for 
 problem \eqref{tf1}-\eqref{tf2} if the following conditions are satisfied:
 \par
 1) $\mathcal A^{tr}$ is compact in $\Theta_+^{weak}$;
 \par
 2) It is strictly invariant with respect to the trajectory semigroup: $T(h)\mathcal A^{tr}=\mathcal A^{tr}$;
 \par
 3) It attracts the images of all bounded sets of $\mathcal K_+$ as time tends to infinity, i.e., 
 for every bounded subset $B\subset\mathcal K_+$ and every neighborhood $\mathcal O(\mathcal A^{tr})$ 
 (in the topology of $\Theta_+^{weak}$), there exists a time $T=T(B,\mathcal O)$ such that
 $$
   T(h)B\subset\mathcal O(\mathcal A^{tr}),
 $$
 for all $h\ge T$.
\edde
Next, we can state the existence result for the above introduced object.
\bete\label{ThW.attr} 
 Let the assumptions of\/ {\rm Theorem \ref{teodiss}} hold. Then, problem\/
 \eqref{tf1}-\eqref{tf2} possesses a trajectory attractor $\mathcal A^{tr}$ in 
 the sense of the above definition. Moreover, this attractor is generated by 
 all complete (i.e., defined for all $t\in\R$) and bounded trajectories for that system.
 Namely, we have
 \begin{equation}\label{W.trdesc}
   \mathcal A^{tr}:=\mathcal K\big|_{t\ge0},
 \end{equation}
 where $\mathcal K\subset L^\infty(\R,\mathcal E)$ is the set of all solutions of 
 \eqref{tf1}-\eqref{tf2} which satisfy
 $$
   \mathcal E(u(T))
    + \int_T^{T+1}\delta\|u_t(t)\|^2
    + \|b^{1/2}(u(t))\Nx w(t)\|^2
    + \|\Dx u(t)\|^2
    + (f'(u(t))\Nx u(t),\Nx u(t))
     \,\dit \le C_*,
 $$
 for all $T\in\R$ and some $C_*>0$.
\ente
\begin{proof} 
As usual (see \cite{CV}), in order to show the attractor existence, we only need 
to verify the existence of a compact and bounded absorbing set for the trajectory 
dynamical system $T(h):\mathcal K_+\to\mathcal K_+$ (the continuity of the semigroup 
in the $\Theta_+^{weak}$-topology is obvious since it is just a translation semigroup). Note that, due to
\eqref{W.en1}, the set $\mathcal B\subset \mathcal K_+$ of solutions $u$ satisfying 
\begin{align}\no
  & \mathcal E(u(T))
   + \int_T^{T+1}\delta\|u_t(t)\|^2
   + \|b^{1/2}(u(t))\Nx w(t)\|^2\\
 \label{W.en2}
  & \mbox{}~~~~~~~~~~~~~~~~~~~~~~~~~
   + \|\Dx u(t)\|^2
   + (f'(u(t))\Nx u(t),\Nx u(t))\,\dit \le 2C_*,
\end{align}
for all $T\ge0$ will be an absorbing set for the semigroup $T(h)$ acting on $\mathcal K_+$. 
Obviously, this set is bounded (in the sense of the Definition~\ref{DefW.bound}). Thus, we only need 
to verify that it is compact in the $\Theta_+^{weak}$ topology.

Indeed, let $\{u_n\}\subset \mathcal B$ be a sequence of solutions. Then, due to estimate \eqref{W.en2}, 
this sequence is precompact in $\Theta_+^{weak}$, so, without loss of generality, we may assume 
that $u_n\to u\in \Theta_+^{weak}$ in the topology of $\Theta_+^{weak}$ and we only need to verify 
that the limit function $u$ solves \eqref{tf1}-\eqref{tf2} and satisfies \eqref{W.en2} as well.

The proof of this fact repeats almost word by word the proof of the existence Theorem \ref{teoesi}
and is even a bit simpler since we do not need to consider the regular approximations to $f$ and $b$ 
(note that the uniform estimate \eqref{W.en2} allows us to control uniformly all of the norms 
involved into \eqref{regou}-\eqref{regofu}). By this reason, we leave the rigorous proof to the reader.

Thus, all of the assumptions of the abstract attractor existence theorem 
are verified and the theorem is proved.
\end{proof}


\section{Energy equalities and strong attraction} 
\label{sec:strongattra}

In the viscous case $\delta>0$ and under slightly more
restrictive assumptions on the growth of $f$, we can
prove that $\calA_{\tr}$ is in fact a {\sl strong}\/
trajectory attractor (i.e., it attracts
with respect of the strong topology of $\calX$).
This is the object of our next result:
\bete\label{teointerm}
 Let assumptions\/ \eqref{hpb}-\eqref{hpdelta}
 hold and let, additionally, $\delta>0$
 and $\beta=0$. In addition, let
 \begin{equation}\label{9bis}
   \kappa\ge \frac{3s}2+1~~\text{if }\,d=3,
    \qquext{and~~}\kappa \ge s+1~~\text{if }\,d=2.
 \end{equation}
 Then, {\rm any} weak solution 
 of problem \eqref{tf1}-\eqref{tf2}
 satisfies the additional regularity properties
 \begin{equation}\label{regofunew}
   w\in L^2(0,T;H), \qquad
   f(u)\in L^2(0,T;H).
 \end{equation}
 Moreover, a.e.~in $(0,\infty)$, there holds
 the following energy\/ {\rm equality}:
 \begin{equation}\label{energyeq}
   \ddt\calE(u)
    +\delta \| u_t \|^2
    +\io b(u)|\nabla w|^2
   = 0,
 \end{equation}
 as well as the following entropy\/
 {\rm equality}
 (compare with \eqref{st21}):
 \begin{equation}\label{S.ent}
  \ddt \Big(\io M(u)
   +\frac\delta2\|\nabla u\|^2\Big)
   +\|\Delta u\|^2
   +\io \big(f'(u)+\gamma'(u)\big)|\nabla u|^2
   +\big(g,\Delta u\big)
  =0.
 \end{equation}
\ente
\begin{proof}
Let us start proving \eqref{energyeq} and first
deal with the 3D-case.
The key step is given by the following integration
by parts formula.
\bele\label{lemipepa}
 Let $b\in L^p(\Omega)$ for some $p>1$,
 with $b\ge 0$ a.e.~in~$\Omega$. Let also
 $b^{-1}\in L^q(\Omega)$ for some $q>3/2$
 if $d=3$ (respectively, for some $q>1$
 if $d=2$). Let $\phi\in H$ and let
 $w$ be the (unique) solution to the
 degenerate elliptic problem
 \begin{equation}\label{degelli}
   -\dive(b\nabla w) + w = \phi, \quext{in }\,\Omega,
     \qquad (b\nabla w)\cdot \bn = 0, \quext{on }\,\Gamma.
 \end{equation}
 Then, $b|\nabla w|^2\in L^1(\Omega)$ and
 \begin{equation}\label{ipepa}
   \big( -\dive(b\nabla w), w \big) = \io b |\nabla w|^2.
 \end{equation}
\enle
\beos\label{sulletracce}
 As it will be clear from the proof, the regularity of $w$
 is sufficient to state \eqref{degelli} in that ``strong''
 form. In particular, since we have
 \begin{equation}\label{Lpdiv}
   b\nabla w \in L^{\frac{2p}{p+1}}(\Omega), \qquad
    \dive (b\nabla w)\in H,
 \end{equation}
 a suitable trace theorem (cf., e.g.,
 \cite[Thm.~2.7.6]{BrGi})
 permits to interpret the boundary
 condition in \eqref{degelli}
 in the sense of trace operators as a relation
 in the space $W^{-\frac{p+1}{2p},\frac{2p}{p+1}}(\Gamma)$.
\eddos
\noindent%
{\bf Proof of Lemma~\ref{lemipepa}.~~}%
 Again, we prove the theorem for $d=3$ and just point out some
 minor differences occurring for $d=2$.
 Let $A$ be the Laplace operator with $0$-Neumann boundary conditions,
 namely,
 \begin{equation}\label{defiA}
   A: V\to V', \qquad
    \duav{Av,z}:=\io \nabla v\cdot \nabla z.
 \end{equation}
 Then, one can see $A+\Id$ as a strictly positive unbounded
 operator on $H$ and consider fractional powers of it.
 For $\epsi>0$, we let $b\ee:=\max\{b,\epsi\}$. We now consider
 the approximate problem
 \begin{equation}\label{degelliee}
   \epsi A^3 w\ee - \dive(b\ee\nabla w\ee) + w\ee = \phi,
   \qquad (b\ee\nabla w\ee)\cdot \bn = 0 \quext{on }\,\Gamma.
 \end{equation}
 Then, testing \eqref{degelliee} by $w\ee$
 we obtain
 \begin{equation}\label{elliee15}
   \epsi \| A^{3/2} w\ee \|^2
    + \io b\ee |\nabla w\ee|^2
    + \| w\ee \|^2
    = ( \phi , w\ee ).
 \end{equation}
 Thus, for all $\epsi>0$, we have that
 $w\ee\in D((A+\Id)^{3/2})\subset H^3(\Omega)$.

 From \eqref{elliee15}, we obtain that
 $w\ee$ is bounded, independently of $\epsi$,
 in $H$. In addition, we have
 \begin{equation}\label{elliee11}
   \big\| \nabla w\ee \big\|_{L^{\frac{2q}{q+1}}(\Omega)}
    \le \big\| b\ee^{1/2} \nabla w\ee \big\|
     \big\| b\ee^{-1/2} \big\|_{L^{2q}(\Omega)}
    \le c
 \end{equation}
 and, being $q>3/2$, it follows $2q/(q+1) > 6/5$
 (respectively, if $d=2$, from $q>1$ we have
 $2q/(q+1) > 1$). Thus, by standard compact embedding
 results, we have, up to a
 (nonrelabelled) subsequence of $\epsi\searrow 0$,
 \begin{equation}\label{stee11}
   w\ee \to w \quext{weakly in }\,W^{1,\frac{2q}{q+1}}(\Omega)
    \quext{and strongly in }\,H.
 \end{equation}
 Moreover, being $p>1$, we can write
 \begin{equation}\label{elliee12}
   \| b\ee\nabla w\ee \|_{L^{\frac{2p}{p+1}}(\Omega)}
    \le \| b\ee^{1/2} \nabla w\ee \|
    \| b\ee^{1/2} \|_{L^{2p}(\Omega)}
    \le c.
 \end{equation}
 Thus, $b\ee \nabla w\ee$ is bounded in
 $L^{2p/(p+1)}(\Omega)\subset D((A+\Id)^{-1})$
 and, consequently,
 \begin{equation}\label{elliee14}
   \big\| -\dive( b\ee\nabla w\ee) \big \|_{D((A+\Id)^{-3/2})}
    \le c
 \end{equation}
 and, proceeding similarly with Subsection~\ref{pass-lim},
 we can also prove that $b\ee\nabla w\ee$ tends to $b\nabla w$
 weakly in $L^{\frac{2p}{p+1}}(\Omega)$.
 Moreover, \eqref{elliee14} tells us that,
 for any $\epsi>0$, equation
 \eqref{degelliee} makes sense at least as a relation
 in $D((A+\Id)^{-3/2})$ (in particular, the obtained
 regularity $w\ee\in D((A+\Id)^{3/2})$ justifies having used
 $w\ee$ as a test function in \eqref{degelliee}).

 Now, the obvious fact that $b\ee\to b$ strongly
 in $L^p(\Omega)$, the first of \eqref{stee11},
 and Ioffe's theorem (see, e.g., \cite{ioffe77}) give
 \begin{equation}\label{elliee15bis}
    \io b |\nabla w|^2
     \le \liminf_{\epsi\searrow0} \io b\ee |\nabla w\ee|^2.
 \end{equation}
 Thus, using \eqref{elliee15} and the second of \eqref{stee11},
 we can go on as follows:
 \begin{align}\no
   \io b |\nabla w|^2
    & \le \lim_{\epsi\searrow 0} ( \phi - w\ee , w\ee )
    + \liminf_{\epsi\searrow 0}
       \big( - \epsi \| A^{3/2} w\ee \|^2 \big)\\
    \label{elliee16}
    & \le ( \phi - w , w )
       = \big( -\dive(b\nabla w) , w \big),
 \end{align}
 where \eqref{degelli} has been used to deduce the last equality.
 Thus, to complete the proof, we have to show the inequality
 converse to \eqref{elliee16}.
 Now, $b\nabla w\in L^{2p/(p+1)}(\Omega)$ thanks to
 (the $\liminf$ of) \eqref{elliee12}.
 Thus, also $-\dive(b\nabla w)\in D((A+\Id)^{-3/2})$
 so that we can test \eqref{degelli}
 by $w\ee\in D((A+\Id)^{3/2})$ and
 rigorously integrate by parts to obtain
 \begin{align}\no
   \big( -\dive(b\nabla w) , w\ee \big)
    & = \io b \nabla w \cdot \nabla w\ee
     \le \frac12  \io b |\nabla w|^2
     + \frac12  \io b |\nabla w\ee|^2 \\
   \no
    &  \le \frac12  \io b |\nabla w|^2
     + \frac12  \io b\ee |\nabla w\ee|^2\\
   \label{elliee17}
    &  \le \frac12  \io b |\nabla w|^2
    - \frac{\epsi}2 \| A^{3/2} w\ee \|^2
    + \frac12 ( \phi - w\ee , w\ee ),
 \end{align}
 where the fact that $b\le b\ee$ almost everywhere and
 the equality \eqref{elliee15} have also been used.
 Passing to the limit and using the
 second \eqref{stee11} and that
 $-\dive(b\nabla w)\in H$ (as it follows by
 comparison in \eqref{degelli}),
 we then obtain
 \begin{equation}\label{elliee18}
  \big( -\dive(b\nabla w) , w \big)
   \le \frac12  \io b |\nabla w|^2
     + \frac12 ( \phi - w , w )
   \le \frac12  \io b |\nabla w|^2
     + \frac12 \big ( -\dive(b\nabla w) , w \big).
 \end{equation}
 Namely, we obtained the inequality converse
 to \eqref{elliee16}, whence the thesis.~\dimbox

\medskip

\noindent%
We now proceed with the proof of Theorem~\ref{teointerm}
and, precisely, of equality \eqref{energyeq}
under the assumption \eqref{9bis}. To do this,
we first prove \eqref{regofunew} and,
with this purpose, we set $z:=u^{-1}$ and 
observe that equation \eqref{tf2}
%
%
can be rewritten as
\begin{equation}\label{tfz}
  \delta z_t + z^2 \Delta z^{-1} + z^{\kappa+2}
   = - z^2 w + z^2 \phi,
   \qquext{where }\,\phi:= \gamma(u) - g
\end{equation}
and we notice that
\begin{equation}\label{regophi}
  \| \phi \|_{L^\infty(\Omega\times(0,T))}
   \le C,
\end{equation}
thanks to \eqref{hpgamma}-\eqref{hpg}. Here and
below, $C$ denotes a constant possibly
depending on the ``energy'' of the initial data
(cf.~\eqref{hpu0}) and on the choice of $T$, while
$c$ is an absolute constant (i.e., it does not
depend on the initial data or on $T$).
By \eqref{regosergey}, we also know that
\begin{equation}\label{regozw}
  \| u^{s/2} \nabla w \|_{L^2(\Omega\times(0,T))}
   \le C.
\end{equation}
%
%
%
%
%
%
%
Next, by the first of \eqref{regofu},
a comparison in \eqref{tf2} gives also
\begin{equation}\label{stwl2l1}
  \| w \|_{L^2(0,T;L^1(\Omega))}
   \le C.
\end{equation}
At this point, we note that, for $d=3$,
thanks to \eqref{regoenergy} and \eqref{9bis},
\begin{equation}\label{16lug11}
  \| b^{-1/2}(u) \|_{L^\infty(0,T;L^3(\Omega))}
   \le C.
\end{equation}
Hence, combining \eqref{16lug11} with \eqref{stwl2l1}
and \eqref{regosergey}, we readily arrive at
\begin{equation}\label{16lug12}
  \| w \|_{L^2(Q)}
  \le c \| w \|_{L^2(0,T;W^{1,6/5}(\Omega))}
   \le C,
\end{equation}
whence \eqref{regofunew} follows simply by comparing
terms in \eqref{tf2} and taking advantage
of \eqref{regou}-\eqref{regofu}.
In the case $d=2$, we only have the
$L^\infty(0,T;H)$-norm in \eqref{16lug11},
but \eqref{regofunew} still follows
by using the continuous embedding $W^{1,1}(\Omega)\subset H$
in the analogue of \eqref{16lug12}.
%
%
%
%
%
%
%
%
%
%
%
%
%
%
%
%

\smallskip

We now proceed with the proof of \eqref{energyeq}.
Summing together \eqref{tf1} and \eqref{tf2}, we have
\begin{equation}\label{tf1-2}
  - \dive (b(u)\nabla w) + w = \phi:=
   (\delta-1) u_t - \Delta u + f(u) + \gamma(u) - g
\end{equation}
and it is clear from \eqref{hpgamma}-\eqref{hpg},
\eqref{regou} and \eqref{regofunew}
that, for any $T>0$, it is
$\|\phi\|_{L^2(0,T;H)}\le C_T$.

Moreover, from \eqref{hpb2} and (the limit of)
\eqref{st31} (or \eqref{st31-2d}),
we have that, a.e.~in~$(0,T)$,
$b(u)\in L^p(\Omega)$ for a suitable $p>1$
(e.g., if $d=3$, we can take $p=10/s$).
Finally, thanks to \eqref{regofunew}, it
is clear that (if $d=3$, the case $d=2$ 
being analogous),
a.e.~in~$(0,T)$,
\begin{equation}\label{elliee21}
  \| b^{-1}(u) \|_{L^q(\Omega)}
   \le \| z^s \|_{L^q(\Omega)}
   \le c \| z^{\frac{2(\kappa-1)}3} \|_{L^q(\Omega)} + c
    \le C,
\end{equation}
for a suitable $q>3/2$. 
Hence, $b=b(u)$ and $\phi$ 
defined in \eqref{tf1-2}
satisfy, a.e.~in~$(0,T)$, the assumptions
of Lemma~\ref{lemipepa}. Consequently,
as we test \eqref{tf1} by $w$ and \eqref{tf2} by $u_t$,
we obtain, thanks to \eqref{ipepa},
the energy equality \eqref{energyeq}, as desired.

\medskip

Finally, we come to the proof of \eqref{S.ent}.
From \eqref{regofu} and \eqref{9bis}, we see that at least $m(u)\in L^2(Q)$,
so since due to \eqref{regou} $u_t\in L^2(Q)$ as well, we have
\begin{equation}\label{S.1}
  \ddt (M(u),1)
   = (u_t, m(u))
   = \big(\dive (b(u)\nabla w), m(u)\big),
\end{equation}
where the right-hand side is understood as a scalar product in $L^2(Q)$. Thus, we  need to verify that
\begin{equation}\label{S.parts}
  \big(\dive (b(u)\nabla w), m(u)\big)
   = - \big(b(u)\nabla w,\nabla m(u)\big)
   = - (\nabla w,\nabla u)
   = (w,\Dx u),
\end{equation}
{\sl almost}\/ everywhere in time. To this end, we note that, 
according to \eqref{regou} and the embedding $H^2(\Omega)\subset C(\barO)$, we have
$u(t),b(u(t))\in C(\barO)$ for almost all $t$. Then, keeping in mind 
that $b^{1/2}(u)\nabla w\in L^2(Q)$ by \eqref{regosergey},
we conclude that $b(u(t))\nabla w(t)\in L^2(\Omega)$ 
for almost all $t$. 

Thus, we only need to check that $m(u(t))\in H^1(\Omega)$ for almost all $t$. 
The fact that this function belongs to $L^2(\Omega)$ is already verified, so we need that 
$$
  \nabla m(u(t))=b^{-1}(u(t))\nabla u(t)\in L^2(\Omega). 
$$
Since, due to \eqref{regou}, we know that $\nabla u(t)\in L^6(\Omega)$ for 
almost all $t$, it is sufficient to check that $b^{-1}(u(t))\in L^3(\Omega)$. 
Actually, this follows immediately from the proved fact that $f(u)\in L^2(Q)$
(see the proof of the energy equality) and condition \eqref{9bis}.

Thus, we have verified that, for almost all $t$, $m(u(t))\in H^1(\Omega)$ 
and $b(u)\nabla w(t)\in L^2(\Omega)$. This justifies the first two 
equalities in \eqref{S.parts}. Note that the last one is obvious since 
$w(t)\in L^2(\Omega)$ and $\Dx u(t)\in L^2(\Omega)$ for almost all $t$. 
Thus, we have verified that
$$
  \ddt (M(u(t)),1)=(w(t),\Dx u(t))
$$
for almost all $t$. Inserting the expression for $w(t)$ from \eqref{tf2} 
into this identity, we end up with the desired entropy equality \eqref{S.ent}. 
Theorem~\ref{teointerm} is proved.
\end{proof}
%
%
%
%
%
%
%
%
\noindent%
Our next task here is to obtain stronger results on the attraction to the above 
constructed trajectory attractor under the additional assumptions of 
Theorem~\ref{teointerm}. We start by stating a couple of 
corollaries that improve the results of Theorems \ref{teoesi} 
and \ref{teodiss} and simplify the construction of the trajectory phase
space based on the energy and entropy inequalities obtained above.
\beco\label{CorS.dis} 
 Let the assumptions of\/ {\rm Theorem \ref{teointerm}} hold. Then, every weak solution 
 of problem\/ \eqref{tf1}-\eqref{tf2} satisfies the dissipative estimate\/ \eqref{W.en} 
 and, therefore, every weak solution automatically satisfies\/ \eqref{W.en1} with 
 $C_u=Q(\mathcal E(u(0))$. Thus, the condition\/ \eqref{W.en1} in the definition of 
 the trajectory phase space $\mathcal K_+$ can be omitted and we may naturally 
 consider $\mathcal K_+$ just as the set of all weak solutions of problem\/
 \eqref{tf1}-\eqref{tf2}. In addition, for every weak solution $u$, we 
 have $u\in C([0,T],H^1(\Omega))$ and $u^{1-\kappa}\in C([0,T],L^1(\Omega))$.
\enco
\noindent%
Indeed, we only need the entropy and energy (in)equalities in order to derive 
the dissipative estimate \eqref{W.en}. Since these inequalities now hold for 
every weak solution, we have this estimate for every weak solution as well. 
The continuity properties stated in the corollary 
follow immediately from the energy equality.
\beco\label{CorS.bound} 
 Let the assumptions of\/ {\rm Theorem \ref{teointerm}} hold. 
 Then, the set $B$ is bounded 
 in $\mathcal K_+$, in the sense of {\rm Definition \ref{DefW.bound}},
 if and only if the set of the initial data $\{u(0),\ u\in B\}$ 
 is bounded in the energy space $\calX_\mu$.
\enco
\noindent%
We are now able to state our main result on the strong convergence to the trajectory attractor.

\bete\label{ThS.strong} 
 Let the assumptions of\/ {\rm Theorem \ref{teointerm}} hold. 
 Then, the trajectory attractor 
 $\mathcal A^{tr}$ of problem \eqref{tf1}-\eqref{tf2} is compact in 
 $C_{loc}(\R_+,\mathcal X)$ and the attraction property holds in that strong topology
 as well (remind that the compactness in $C_{loc}(\R_+,\mathcal X)$ means that the 
 $u$-component of $\mathcal A^{tr}$ is compact in $C_{loc}(\R_+,H^1(\Omega))$ and 
 the $u^{1-\kappa}$-component is compact in $C_{loc}(\R_+,L^1(\Omega))$).
\ente
\begin{proof} 
Let $\mathcal B\subset K_+$ be the absorbing set introduced 
in the proof of Theorem \ref{ThW.attr}.
We claim that the set 
$$
  \mathcal B_1:=T(1)\mathcal B
$$
is an absorbing set which is compact in the above mentioned topology (this is clearly 
enough for the proof of the theorem). Indeed, let $\{u_n\}\subset\mathcal B$ be an arbitrary 
sequence of solutions. Then, since $\mathcal B$ is compact in $\Theta_+^{weak}$, we 
may assume without loss of generality that $u_n\to u$ in $\Theta_+^{weak}$, where $u$
also solves the problem \eqref{tf1}-\eqref{tf2}. To verify the above mentioned 
compactness, we need to check that
\begin{equation}\label{nuovaconv}
  u_n\to u\ \text{ in }\, C([1,N],H^1(\Omega)),\qquad
   u_n^{1-\kappa}\to u^{1-\kappa} \ \text{ in }\, C([1,N], L^1(\Omega)),
\end{equation}
for every $N>1$. Furthermore, without loss of generality, we
may check these convergences for $N=2$ only.

To this end, we will use the proved energy equality which we will rewrite in the following form:
\begin{align}\no
  & T\Big [ \frac12 \|\nabla u_n(T)\|^2
   + (F(u_n(T)),1)
   + (\Gamma(u_n(T)),1)
   - (g,u_n(T))\Big] \\
 \no
  & \mbox{}~~~~~~~~~~
   + \int_0^T \delta t \|\Dt u_n(t)\|^2
    + t \big(b(u_n(t)) \nabla w_n(t),\nabla w_n(t)\big)\,\dit\\
 \label{S.en}
  & \mbox{}~~~~~
  = \int_0^T \frac12\|\nabla u_n(t)\|^2
   + (F(u_n(t)),1)
   + (\Gamma(u_n(t)),1)
   - (g,u_n(t))
    \,\dit,
\end{align}
where $T\in[1,2]$. Our next task is to pass to the limit $n\nearrow\infty$ in this 
inequality. First of all, thanks to the energy and entropy estimates and to
the Aubin-Lions compactness theorem, we have
\begin{equation}\label{S.-1}
  u_n \to u \quext{in }\,C_w([1,2];V),
   \qquext{(strongly) in }\,C([1,2];H)\cap L^2(1,2;V),
\end{equation}
and pointwise (a.e.). Thus, using the first convergence above and 
Fatou's Lemma, we see that
\begin{equation}\label{S.000}
  \frac12 \|\nabla u(T)\|^2 \le\liminf_{n\nearrow\infty} 
      \frac12 \|\nabla u_n(T)\|^2, \qquad
   (F(u(T)),1)\le\liminf_{n\nearrow\infty} (F(u_n(T)),1).
\end{equation}
Next, thanks also to \eqref{hpgamma}-\eqref{hpg},
\begin{equation}\label{S.new1}
  (\Gamma(u),1)-(g,u)
   = \lim_{n\nearrow\infty} (\Gamma(u_n),1)-(g,u_n)
   \quext{strongly in }\,C^0([1,2]).
\end{equation}
Moreover, thanks to lower semicontinuity of 
norms w.r.t.~weak convergence and to Ioffe's theorem,
we also have
\begin{align}
  & \int_0^T\delta t\|\Dt u(t)\|^2\,\dit
   \le\liminf_{n\nearrow\infty} \int_0^T\delta t\|\Dt u_n(t)\|^2\,\dit,\\ 
  & \int_0^T t \big(b(u(t))\nabla w(t),\nabla w(t)\big)\,\dit
   \le \liminf_{n\nearrow\infty} 
      \int_0^T t \big (b(u_n(t))\nabla w_n(t),\nabla w_n(t)\big)\,\dit.
\end{align}
Finally, in order to pass to the limit in \eqref{S.en}, we only need to prove that 
$F(u_n)\to F(u)$ {\sl strongly} in $L^1([0,2]\times\Omega)$. 
Actually, this follows from the uniform $L^2$-bound of $f(u_n)$,
the pointwise convergence $u_n\to u$ and the generalized Lebesgue theorem.
Thus, we can take the supremum limit $n\nearrow\infty$ in
\eqref{S.en} and obtain the {\sl inequality}:
\begin{align}\no
  & \limsup_{n\nearrow\infty} 
   T\Big [ \frac12 \|\nabla u_n(T)\|^2
   + (F(u_n(T)),1) \Big] \\
  \no
 & \mbox{}~~~~~ \le  
  -T \big[ (\Gamma(u(T)),1)
   - (g,u(T)) \big]
   - \int_0^T \delta t \|\Dt u(t)\|^2
   - t \big(b(u(t)) \nabla w(t),\nabla w(t)\big)\,\dit\\
 \label{S.en2}
  & \mbox{}~~~~~~~~~~
   + \int_0^T \frac12\|\nabla u(t)\|^2
   + (F(u(t)),1)
   + (\Gamma(u(t)),1)
   - (g,u(t))  \,\dit.
\end{align}
On the other hand, applying the energy equality 
in the form~\eqref{S.en}
directly to the limit solution $u$, 
and comparing with \eqref{S.en2}, we obtain, for all~$T\in(1,2)$,
$$ 
  \limsup_{n\nearrow\infty} T\Big [ \frac12 \|\nabla u_n(T)\|^2
   + (F(u_n(T)),1) \Big] 
  \le T\Big [ \frac12 \|\nabla u(T)\|^2
   + (F(u(T)),1) \Big],
$$
whence, recalling \eqref{S.000}, we infer
$$ 
  \|u_n(T)\|_{V} \to \|u(T)\|_{V}, \qquad
    \|F(u_n(T))\|_{L^1(\Omega)} \to \|F(u(T))\|_{L^1(\Omega)}.
$$
This, together with the weak convergence $u_n(T)\to u(T)$ in $V$ 
(cf.~\eqref{S.-1}) and $u_n\to u$ 
almost everywhere, implies the strong convergence
\begin{equation}\label{S.almost}
  u_n(T)\to u(T)\ \ \text{in }\, V, \qquad
   u_n^{1-\kappa}(T)\to u^{1-\kappa}(T)\ \ \text{in\ }\, L^1(\Omega), 
\end{equation}
for all $T\in[1,2]$. This gives the strong convergence
$u_n\to u$ in ${\mathcal X}$ pointwise in time.
The desired uniform convergence \eqref{nuovaconv}
can be easily obtained using the standard contradiction 
arguments and applying the energy equality for $u_n(T_n)$ instead of $u_n(T)$. 
Theorem \ref{ThS.strong} is proved. 
\end{proof}
\beco\label{CorS.last} 
 Arguing in a similar way (and using also the entropy equality), one can verify 
 the compactness and strong convergence to the trajectory attractor in all 
 spaces involved in \eqref{regou}-\eqref{regofu}.
\enco


\section{Separation from singularities and uniqueness}
\label{sec-long}

In this section we prove that, 
in the viscous case $\delta>0$,
if $\kappa$ is large enough, then any weak solution 
becomes uniformly
strictly positive for any $t>0$.
This is the object of the following
\bete\label{teosep}
 Let assumptions\/ \eqref{hpb}-\eqref{hpdelta}
 hold and let, additionally, $\delta>0$
 and $\beta=0$. In addition, let
 \begin{equation}\label{9}
   \kappa>2s+3~\,\text{if }\,d=3,
    \qquext{and~}\,\kappa>s+1\ge 2~~\text{if }\,d=2.
 \end{equation}
 Then, there exists a function
 $Q;[0,\infty)^2\to [0,\infty)$, monotone in
 each of its arguments, such that any weak solution $u$ 
 satisfies, for any $\epsilon>0$, the\/ 
 {\rm separation property}
 \begin{equation}\label{10}
   \|u^{-1}(t)\|_{L^\infty(\Omega)}\le
    Q\big(\EE_0,\epsilon^{-1}\big)
    \quext{for a.e.~}\/ t\ge \epsilon.
 \end{equation}
\ente
\noindent%
The proof of the theorem will be given later 
in this section. As a consequence, we also have
further time-regularization properties that
imply uniqueness for strictly positive times
as well:
\bete\label{teouni}
 Let the assumptions of\/ {\rm Theorem~\ref{teosep}}
 hold (in particular, let $\delta>0$). Then, 
 for any $\epsilon>0$ and any weak solution $u$ there holds:
 \begin{align}\label{contouni12}
   & w\in L^2(\epsilon,T;V),\\
  \label{contouni11}
   & u\in H^1(\epsilon,T;V)\cap L^\infty(\epsilon,\infty;H^2(\Omega)).
 \end{align}
 Moreover, in the class of weak solutions uniqueness holds
 at least for strictly positive times.
\ente
\beos\label{arbireg}
 We point out that \eqref{contouni12}-\eqref{contouni11},
 which {\sl suffice}\/ to prove uniqueness, are however not presumed to be
 optimal properties. Actually, thanks to \eqref{10},
 \eqref{tf1} is nondegenerate and \eqref{tf2} is nonsingular
 for strictly positive times. Thus, by means of classical methods, 
 one could easily prove that the solution $u$ becomes,
 instantaneously in time, arbitrarily regular, provided of course
 that also the data $\gamma$ and $g$ are smooth.
\eddos
\noindent%
As a consequence of uniqueness, we finally
have
\beco\label{teoattr}
 Let the assumptions of\/ {\rm Theorem~\ref{teosep}}
 hold (in particular, let $\delta>0$). Then, the 
 dynamical process generated by weak solutions
 admits the (strong) global attractor $\calA$
 in the standard sense (to be more precise,
 in the sense of semigroups with unique continuation). Namely,
 $\calA$ is a compact and fully invariant subset of $\calX_\mu$
 such that, for any bounded set $B\subset \calX_\mu$ there
 holds
 \begin{equation}\label{attra}
   \lim_{t\nearrow\infty}
    \deriv_{\calX}(u(t),\calA) = 0,
 \end{equation}
 uniformly with respect to weak solutions $u$ such
 that $u(0)\in B$.
\enco
%


\subsection{Proof of Theorem~\ref{teosep} in the 3D-case}
\label{proof-sep}

We consider equation \eqref{tf2} rewritten in the form
\eqref{tfz}, where we assume $\delta=1$ for simplicity.
We also assume $s>0$, the case $s=0$ being simpler
since one can directly take advantage of the 
$L^2(0,T;L^6(\Omega))$-regularity of $w$.
Then, the proof is based on a suitable
version of the Moser iteration argument, i.e., we will
take $\nu>1$ and test \eqref{tfz} by $\nu z^{\nu-1}$
for increasing exponents $\nu$.
We have to remark that this procedure, apparently having 
a formal character since the above test function 
could grow very fast and, hence, have insufficient
regularity, can be easily justified simply by truncating
$z$ at some level $K$ and then letting $K\nearrow\infty$.
In particular, the argument does not require any approximation
of the equation 
and, hence, works for {\sl all}\/ weak solutions
in the class introduced in Definition~\ref{def:weaksol}

That said, testing \eqref{tfz} by $\nu z^{\nu-1}$ and
integrating over $Q_\nu:=\Omega\times(\tau_\nu,T)$,
where the ``initial'' time
$\tau_\nu$ will be chosen later on, we then have
\begin{equation}\label{contonu11}
  J_\nu^\nu
   + \iint_{Q_\nu} z^{\kappa+\nu+1}
  \le \| z(\tau_\nu) \|_{L^\nu(\Omega)}^\nu
   + \nu \iint_{Q_\nu} \phi z^{\nu+1}
   - \nu \iint_{Q_\nu} w z^{\nu+1},
\end{equation}
where we have set
\begin{equation}\label{defiJnu}
  J_\nu^\nu
   := \| z \|_{L^\infty(\tau_\nu,T;L^\nu(\Omega))}^\nu
  + \| \nabla z^{\nu/2} \|_{L^2(\tau_\nu,T;H)}^2.
\end{equation}
Adding now
\begin{equation}\label{contonu12}
   \| z^{\nu/2} \|_{L^2(\tau_\nu,T;H)}^2
    = \| z \|_{L^{\nu}(Q_\nu)}^\nu
\end{equation}
to both hands sides of \eqref{contonu11},
in order to recover the full $V$-norm
of $z^{\nu/2}$, and setting
\begin{equation}\label{defiInu}
  I_\nu^\nu:=
   \| z \|_{L^\infty(\tau_\nu,T;L^\nu(\Omega))}^\nu
   + c\OO \| z \|_{L^\nu(\tau_\nu,T;L^{3\nu}(\Omega))}^\nu
   \ge \| z \|_{L^{5\nu/3}(Q_\nu)}^\nu,
\end{equation}
where $c\OO$ is a suitable embedding constant,
we then arrive at
\begin{equation}\label{contonu13}
  I_\nu^\nu
   + \iint_{Q_\nu} z^{\kappa+\nu+1}
  \le \| z(\tau_\nu) \|_{L^\nu(\Omega)}^\nu
   + \nu C \| z \|_{L^{\nu+1}(Q_\nu)}^{\nu+1}
   + \| z \|_{L^{\nu}(Q_\nu)}^\nu
   - \nu \iint_{Q_\nu} w z^{\nu+1},
\end{equation}
and we have to provide a bound for the last term
on the \rhs.

To do this, let $(p,p^*)$ and $(q,q^*)$ be
two couples of conjugate exponents with $p\le 2$. Then,
using also \eqref{stwl2l1}, we obtain
\begin{align}\no
  - \nu \iint_{Q_\nu} w z^{\nu+1}
   & \le \nu \| w \|_{L^p(\tau_\nu,T;L^{q}(\Omega))}
    \| z^{\nu+1} \|_{L^{p^*}(\tau_\nu,T;L^{q^*}(\Omega))} \\
  \no
   & \le \nu \big( \| w \|_{L^p(\tau_\nu,T;L^1(\Omega))}
    + \| \nabla w \|_{L^p(Q_\nu)} \big)
    \| z^{\nu+1} \|_{L^{p^*}(\tau_\nu,T;L^{q^*}(\Omega))}\\[1mm]
  \label{contonu14}
   & \le \nu \big( C
    + \| \nabla w \|_{L^p(Q_\nu)} \big)
    \| z^{\nu+1} \|_{L^{p^*}(\tau_\nu,T;L^{q^*}(\Omega))},
\end{align}
provided that we choose $q=3p/(3-p)$,
so that $W^{1,p}(\Omega)\subset L^q(\Omega)$ continuously.

Then, the term with $\nabla w$ is estimated
this way:
\begin{align}\no
  \| \nabla w \|_{L^p(Q_\nu)}
   & = \| u^{s/2} z^{s/2} \nabla w \|_{L^p(Q_\nu)}
    \le \| u^{s/2} \nabla w \|_{L^2(Q_\nu)}
    \| z^{s/2} \|_{L^{\frac{2p}{2-p}}(Q_\nu)}\\
 \label{contonu15}
   & \le C \| z^{s/2} \|_{L^{\frac{2p}{2-p}}(Q_\nu)},
\end{align}
thanks also to \eqref{regozw}. Thus,
collecting \eqref{contonu14} and \eqref{contonu15},
we have
\begin{equation}\label{contonu16}
  - \nu \iint_{Q_\nu} w z^{\nu+1}
   \le C \nu \Big( 1 + \| z^{s/2} \|_{L^{\frac{2p}{2-p}}(Q_\nu)} \Big)
   \| z \|_{L^{p^*(\nu+1)}(\tau_\nu,T;L^{q^*(\nu+1)}(\Omega))}^{\nu+1}.
\end{equation}
Now, let us assume to know a bound of the term
$I_{n-1}=I_{\nu_{n-1}}$ from the preceding step of the
iteration. Then, thanks to the last inequality
in \eqref{defiInu}, we can use it to estimate the term
in brackets so to have
\begin{equation}\label{contonu17}
  - \nu \iint_{Q_\nu} w z^{\nu+1}
   \le C \nu \big( 1 + I_{n-1}^{s/2} \big)
   \| z \|_{L^{p^*(\nu+1)}(\tau_\nu,T;L^{q^*(\nu+1)}(\Omega))}^{\nu+1},
\end{equation}
provided that one chooses $p$ as follows:
\begin{equation}\label{contonu18}
  \frac{p}{2-p} = \frac{5\nu_{n-1}}{3s},
   \quext{i.e., }\,
   \frac1p = \frac12 + \frac{3s}{10\nu_{n-1}}.
\end{equation}
This gives in turn
\begin{equation}\label{contonu19}
  \frac1q = \frac16 + \frac{3s}{10\nu_{n-1}}, \qquad
  \frac1{p^*} = \frac12 - \frac{3s}{10\nu_{n-1}}, \qquad
  \frac1{q^*} = \frac56 - \frac{3s}{10\nu_{n-1}}.
\end{equation}
Then, we have to take $\nu=\nu_n$ in a way suitable
for the next step of the iteration. The choice is dictated
by the exponents of the last term in~\eqref{contonu17};
namely, $\nu=\nu_n$ should be close enough to $\nu_{n-1}$
in order that term be still controlled by $I_{n-1}$.
Using interpolation, we require that, for some
$\theta\in [0,1]$,
\begin{equation}\label{contonu20}
  \frac1{p^*(\nu_n+1)}
    =\frac{1-\theta}{\infty}+\frac\theta{\nu_{n-1}},
    \qquad
  \frac1{q^*(\nu_n+1)}
    =\frac{1-\theta}{\nu_{n-1}}+\frac\theta{3\nu_{n-1}}.
\end{equation}
To compute $\theta$, we first take the quotient
of the above equalities and then use \eqref{contonu19}.
This gives
\begin{equation}\label{contonu21}
  \frac{3-2\theta}{3\theta}
   =\frac{p^*}{q^*}
   =\frac{25\nu_{n-1}-9s}{15\nu_{n-1}-9s},
\end{equation}
whence
\begin{equation}\label{contonu22}
  \theta=\frac{15\nu_{n-1}-9s}{35\nu_{n-1}-15s}, \qquand
    \nu_n=\frac{\nu_{n-1}}{p^*\theta}-1
     =\frac{7\nu_{n-1}-3s-6}6,
\end{equation}
where the second of \eqref{contonu19} has also been used.

Thus, it turns out that $\nu_n>\nu_{n-1}$ provided
that $\nu_{n-1}>3(s+2)$. Thus, in order the above iteration
could be performed, we need to find some $\barn\in\NN$
and some $\nu_{\barn}>3(s+2)$
such that, for any $\epsilon\in(0,1)$, there holds
\begin{equation}\label{regobarn}
  I_{\barn} = \Big(  \| z \|_{L^\infty(\epsilon,T;L^{\nu_{\barn}}(\Omega))}^{\nu_{\barn}}
   + c\OO \| z \|_{L^{\nu_{\barn}}(\epsilon,T;L^{3\nu_{\barn}}(\Omega))}^{\nu_{\barn}} \Big)
         ^{\frac{1}{\nu_{\barn}}}
   \le Q(\epsilon^{-1}),
\end{equation}
where $Q$ is a computable monotone function
(whose expression can depend on the magnitude of the
initial data and of $T$).

Let us pospone the verification of \eqref{regobarn}
and let us now see that, for $n > \barn$, the
induction principle can be applied. Coming back to
\eqref{contonu13}, we then have
\begin{equation}\label{contonu23}
  I_n^{\nu_n}
   \le \| z(\tau_n) \|_{L^{\nu_n}(\Omega)}^{\nu_n}
   + C \nu_n \big( 1 + \| z \|_{L^{\nu_n+1}(Q_n)}^{\nu_n+1} \big)
   + C \nu_n I_{n-1}^{\frac{s}2+\nu_n+1},
\end{equation}
where we wrote $n$ in place of $\nu_n$ in some subscripts
and assumed w.l.o.g.~$I_{n-1}\ge 1$.
Thus, extracting the $\nu_n$-th root
and noting that $\nu_n+1\le 5\nu_{n-1}/3$,
we obtain
\begin{equation}\label{contonu24}
  I_n
   \le \| z(\tau_n) \|_{L^{\nu_n}(\Omega)}
   + (C\nu_n)^{\frac1{\nu_n}} I_{n-1}^{\eta_n},
    \quext{where }\,\eta_n:=\frac{2\nu_n+s+2}{2\nu_n}
\end{equation}
and $C$ is independent of $n$. Moreover,
for (arbitrarily small) $\epsilon\in(0,1)$,
given $\tau_{n-1}$ we can choose
$\tau_n\in [\tau_{n-1},\tau_{n-1}+\epsilon n^{-2}\big]$
such that
\begin{equation}\label{contonu24b}
   \| z(\tau_n) \|_{L^{\nu_n}(\Omega)}^{\frac{5\nu_{n-1}}3}
   \le c \| z(\tau_n) \|_{L^{\frac{5\nu_{n-1}}3}(\Omega)}^{\frac{5\nu_{n-1}}3}
    \le c\frac{n^2}\epsilon \int_{\tau_{n-1}}^{\tau_{n-1}+\epsilon n^{-2}}
     \| z \|_{L^{\frac{5\nu_{n-1}}3}(\Omega)}^{\frac{5\nu_{n-1}}3}
     \le c\frac{n^2}\epsilon I_{n-1}^{\frac{5\nu_{n-1}}3}.
\end{equation}
Thus, \eqref{contonu24} can be rewritten as
\begin{equation}\label{contonu24c}
  I_n \le \Big[ \Big( c \frac{n^2}{\epsilon} \Big)^{\frac{3}{5\nu_{n-1}}}
  + (C\nu_n)^{\frac1{\nu_n}} \Big]
   I_{n-1}^{\eta_n},
\end{equation}
whence a standard computation permits to pass to the limit
w.r.t.~$n\nearrow\infty$. Since $\lim_{n\nearrow\infty} \tau_n$
exists and is less or equal than $c\epsilon$, we then obtain
\begin{equation}\label{contonu25}
  \| u \|_{L^\infty(\Omega\times(\epsilon,T))}
   \le Q(\epsilon^{-1}),
\end{equation}
for $Q$ as in \eqref{regobarn}, as desired.

Thus, to conclude the proof it only remains
to check that \eqref{regobarn} holds.
To do this, we come back to \eqref{contonu13} and use now
the $L^{\kappa+\nu+1}$-norm to estimate the \rhs.
Proceeding as above, we still arrive at \eqref{contonu16}, where
now we have to take
\begin{equation}\label{contonu18new}
  \frac{p}{2-p} = \frac{\kappa + \nu_{n-1} + 1}s
   \quext{i.e., }\,
   \frac1p = \frac12 \Big( 1 + \frac{s}{\kappa + \nu_{n-1} + 1} \Big).
\end{equation}
Thus, we obtain
\begin{equation}\label{contonu19new}
  \frac1q = \frac12 \Big( \frac13 + \frac{s}{\kappa + \nu_{n-1} + 1} \Big), \qquad
  \frac1{p^*} = \frac12 \Big( 1 - \frac{s}{\kappa + \nu_{n-1} + 1} \Big), \qquad
  \frac1{q^*} = \frac12 \Big( \frac53 - \frac{s}{\kappa + \nu_{n-1} + 1} \Big),
\end{equation}
whence we get the analogue of \eqref{contonu17}, i.e.,
\begin{equation}\label{contonu17new}
  - \nu \iint_{Q_\nu} w z^{\nu+1}
   \le C \nu \big( 1 + \Lambda_{n-1}^{s/2} \big)
   \| z \|_{L^{p^*(\nu+1)}(\tau_\nu,T;L^{q^*(\nu+1)}(\Omega))}^{\nu+1},
\end{equation}
where
\begin{equation}\label{defiLanu}
  \Lambda_\nu^\nu
   := \| z \|_{L^\infty(\tau_\nu,T;L^\nu(\Omega))}^\nu
  + \| \nabla z^{\nu/2} \|_{L^2(\tau_\nu,T;H)}^2
  + \| z \|_{L^{\nu + \kappa + 1}(Q_\nu)}^{\nu + \kappa + 1}
\end{equation}
and $\Lambda_n:=\Lambda_{\nu_n}$, as before.
Then, we still have to choose $\nu=\nu_n$ in a suitable way.
Similarly as before, we require that for some $\theta\in [0,1]$ it is
\begin{equation}\label{contonu20new}
  \frac1{p^*(\nu_n+1)}
    =\frac{1-\theta}{\infty}+\frac\theta{\kappa + \nu_{n-1} +1},
    \qquad
  \frac1{q^*(\nu_n+1)}
    =\frac{1-\theta}{\nu_{n-1}} + \frac\theta{\kappa + \nu_{n-1} + 1}.
\end{equation}
To compute $\theta$, we first take the quotient
of the above equalities and then use \eqref{contonu19new}.
This gives
\begin{equation}\label{contonu21new}
  \frac{\theta\nu_{n-1}+(1-\theta)(\nu_{n-1}+\kappa+1)}{\theta\nu_{n-1}}
   =\frac{p^*}{q^*}
   =\frac{5(\nu_{n-1}+\kappa+1)-3s}{3(\nu_{n-1}+\kappa+1-s)},
\end{equation}
whence
\begin{equation}\label{contonu22new}
  \frac1\theta=1+\frac{2\nu_{n-1}}{3(\nu_{n-1}+\kappa+1)-3s}
\end{equation}
and, from the first of \eqref{contonu20new},
\begin{equation}\label{contonu22new2}
  \nu_n=\frac{\nu_{n-1}+\kappa+1}{\theta p^*}-1
   = \frac56\nu_{n-1}+\frac12(\kappa-s-1),
\end{equation}
whence it is clear that $\nu_n>\nu_{n-1}$
if and only if $\nu_{n-1}<3(\kappa-s-1)$.
Then, proceeding similarly
with the previous part of the iteration,
if we start knowing a bound of $\Lambda_{\nu_0}$
for some $\nu_0>1$, then we can reach,
in a finite number $\barn$ of steps,
any $\nu_{\barn}<3(\kappa-s-1)$.
Since we also need $\nu_{\barn}>3(s+2)$
from before,
this leads to the compatibility condition
$3(s+2)<3(\kappa-s-1)$, that is equivalent
to assumption \eqref{9}.

Thus, the proof is concluded provided that
we find $\nu_0>1$ to start the argument.
Actually, we can test \eqref{tfz} by $z^{\iota}$
for small $\iota>0$. We obtain,
for $Q=\Omega\times(0,T)$,
\begin{equation}\label{contonu11new}
  J_{1+\iota}^{1+\iota}
   + \iint_{Q} z^{\kappa+2+\iota}
  \le \| z_0 \|_{L^{1+\iota}(\Omega)}^{1+\iota}
   + \iint_{Q} \phi z^{2+\iota}
   - \iint_{Q} w z^{2+\iota}
\end{equation}
and, being $\kappa>3$ by \eqref{9}, it is clear 
that, at least for $\iota<1$,
\begin{equation}\label{contonu11newb}
  \iint_{Q} \phi z^{2+\iota}
  - \iint_{Q} w z^{2+\iota}
   \le \frac12 \iint_{Q} z^{\kappa+2+\iota}
   + C,
\end{equation}
thanks also to \eqref{regofunew}.
%
%
%
%
%
%
%
%
Hence, we can take $\nu_0=1+\iota$
for arbitrary $\iota\in(0,1)$, which concludes
the proof in the case $d=3$.


\subsection{Proof of Theorem~\ref{teosep} in the 2D-case}
\label{proof-sep2d}

The proof is carried out by the very same scheme used in the
3D case, the differences being limited to the exponents
related to use of interpolation and embeddings. Thus, we limit
ourselves to point out these differences.
Now, in place of \eqref{defiInu}, we have
\begin{equation}\label{defiInu2d}
  I_\nu^\nu:=
   \| z \|_{L^\infty(\tau_\nu,T;L^\nu(\Omega))}^\nu
   + \| z^{\frac{\nu}2} \|_{L^2(\tau_\nu,T;V)}^2
   \ge \| z \|_{L^{2\nu}(Q_\nu)}^\nu.
\end{equation}
Thus, taking $p\in(1,2)$, we have $q=2p/(2-p)$, so that,
to control the \rhs\ of \eqref{contonu17},
we need to choose $p$ so that
\begin{equation}\label{contonu182d}
  \frac{p}{2-p} = \frac{2\nu_{n-1}}s,
   \quext{i.e., }\,
   \frac1p = \frac12 \Big( 1 + \frac{s}{2\nu_{n-1}} \Big),
\end{equation}
whence we obtain
\begin{equation}\label{contonu192d}
  \frac1q = \frac{s}{4\nu_{n-1}}, \qquad
  \frac1{p^*} = \frac12 \Big( 1 - \frac{s}{2\nu_{n-1}} \Big), \qquad
  \frac1{q^*} = 1 - \frac{s}{4\nu_{n-1}}
\end{equation}
and, correspondingly,
\begin{equation}\label{contonu202d}
  \frac1{p^*(\nu_n+1)}
    =\frac{1-\theta}{\infty}+\frac\theta{2 \nu_{n-1}},
    \qquad
  \frac1{q^*(\nu_n+1)}
    =\frac{1-\theta}{\nu_{n-1}}+\frac\theta{2 \nu_{n-1}}.
\end{equation}
Thus,
\begin{equation}\label{contonu212d}
  \frac{2-\theta}{\theta}
   =\frac{p^*}{q^*}
   =\frac{4\nu_{n-1}-s}{2\nu_{n-1}-s},
\end{equation}
whence
\begin{equation}\label{contonu222d}
  \theta=\frac{2\nu_{n-1}-s}{3\nu_{n-1}-s}, \qquand
    \nu_n=\frac32\nu_{n-1}-\frac{s+2}2,
\end{equation}
so that we need to find $\nu_{\barn}>s+2$
in order the procedure works.

To do this, we proceed again as before and,
choosing $p$ as in \eqref{contonu18new}, the other
exponents are then given by
\begin{equation}\label{contonu19new2d}
  \frac1q = \frac{s}{2( \kappa + \nu_{n-1} + 1) }, \qquad
  \frac1{p^*} = \frac12 \Big( 1 - \frac{s}{\kappa + \nu_{n-1} + 1} \Big), \qquad
  \frac1{q^*} = 1 - \frac{s}{2( \kappa + \nu_{n-1} + 1) }.
\end{equation}
Then, taking $\theta\in [0,1]$ as in \eqref{contonu20new},
we now arrive at
\begin{equation}\label{contonu21new2d}
  \frac{\theta\nu_{n-1}+(1-\theta)(\nu_{n-1}+\kappa+1)}{\theta\nu_{n-1}}
   =\frac{p^*}{q^*}
   =\frac{2(\nu_{n-1}+\kappa+1)-s}{\nu_{n-1}+\kappa+1-s},
\end{equation}
whence
\begin{equation}\label{contonu22new2d}
  \frac1\theta=1+\frac{\nu_{n-1}}{\nu_{n-1}+\kappa+1-s},
   \qquand \nu_n = \nu_{n-1}+\frac12(\kappa-1-s),
\end{equation}
so that it is $\nu_n>\nu_{n-1}$ if and only if $\kappa>s+1$,
i.e., \eqref{9} holds. Thus, we can arrive in some 
finite number $\barn$ of steps to have $\nu_{\barn}>s+2$
provided that we can start as before from $\nu_0=1+\iota$
for some (small) $\iota>0$. Actually, we can now take $\iota=\kappa-2$,
which is strictly positive thanks to~\eqref{9}.
Thus, \eqref{contonu11new} can be repeated without any variation
and, of course, we still have \eqref{contonu11newb}
thanks to H\"older's and Young's inequalities.
%
%
%
The proof is complete.


\subsection{Proof of Theorem~\ref{teouni} and Corollary~\ref{teoattr}}
\label{proof-uniq}

Again, we just consider the case $d=3$, the case $d=2$ being simpler.
First of all, we deduce further regularity of weak solutions.
Actually, thanks to \eqref{10}, $u$ is uniformly separated
from $0$ for any time $t\ge\epsilon>0$,
$\epsilon>0$ being arbitrary. Then,
\eqref{tf1} becomes in fact nondegenerate and the energy estimate
gives the improved regularity \eqref{contouni12}.
Moreover, the term $f(u)$ in \eqref{tf2}
is now smooth and we can
apply the linear parabolic theory (or test
\eqref{tf2} by $-(t-\epsilon)\Delta u_t$ and perform standard
computations) to deduce \eqref{contouni11}.

At this point, rewriting \eqref{tf1} as a family of
time-dependent elliptic problems, namely
\begin{equation}\label{tf1ell}
  -\Delta w = \frac1{b(u)} \big(- u_t + b'(u)\nabla u\cdot \nabla w\big),
\end{equation}
relations \eqref{contouni11} and \eqref{contouni12}
permit to see that the \rhs\ belongs to $L^2(\epsilon,T;L^{3/2}(\Omega))$,
whence we obtain
\begin{equation}\label{contouni13}
  w\in L^2(\epsilon,T;W^{2,3/2}(\Omega))
   \subset L^2(\epsilon,T;W^{1,3}(\Omega)).
\end{equation}
To prove uniqueness, we can now consider a couple of
solutions $u_1$, $u_2$, set $u:=u_1-u_2$ (and, correspondingly,
$w:=w_1-w_2$) and take the difference of equations
\eqref{tf1}-\eqref{tf2} to obtain
\begin{align}\label{tf1diff}
  & u_t-\dive (b(u_1)\nabla w) = \dive\big( (b(u_1)-b(u_2)) \nabla w_2\big),\\
 \label{tf2diff}
  & w = \delta u_t - \Delta u + W'(u_1) - W'(u_2),
\end{align}
where $W'=f+\gamma$ can be thought to be globally Lipschitz
in view of the strict positivity of $u_1$ and $u_2$.
Then, we test \eqref{tf1diff} by $w$ and \eqref{tf2diff}
by $u_t$. We obtain, for some $c,\alpha>0$,
\begin{align}\no
  & \frac12 \ddt \| \nabla u \|^2
   + \alpha \| \nabla w \|^2
   + \delta \|u_t\|^2\\
 \label{contouni14}
  & \mbox{}~~~~~
   \le c \| u \| \| u_t \|
  + \io  \big| (b(u_1) - b(u_2)) \nabla w \cdot \nabla w_2 \big|
\end{align}
and we can estimate the last term as follows:
\begin{align}\no
  \io  \big| (b(u_1) - b(u_2)) \nabla w \cdot \nabla w_2 \big|
   & \le \frac\alpha2 \| \nabla w \|^2
    + c \| u \|_{L^6(\Omega)}^2 \| \nabla w_2 \|_{L^3(\Omega)}^2\\
 \label{contouni15}
   & \le \frac\alpha2 \| \nabla w \|^2
    + c \| u \|_{V}^2 \| \nabla w_2 \|_{L^3(\Omega)}^2.
\end{align}
Thanks to \eqref{contouni13}, we can then
apply Gronwall's Lemma to \eqref{contouni14},
which gives the assert. 
At this point, Corollary~\ref{teoattr} is
an immediate consequence of the uniqueness property
and of the general theory of infinite-dimensional dynamical
systems \cite{BVbook,Te}.



\vspace{15mm}

\noindent%
{\bf First author's address:}\\[1mm]
Giulio Schimperna\\
Dipartimento di Matematica, Universit\`a degli Studi di Pavia\\
Via Ferrata, 1,~~I-27100 Pavia,~~Italy\\
E-mail:~~{\tt giusch04@unipv.it}

\vspace{4mm}

\noindent%
{\bf Second author's address:}\\[1mm]
Sergey Zelik\\
Department of Mathematics, University of Surrey\\
Guildford,~~GU2 7XH,~~United Kingdom\\
E-mail:~~{\tt S.Zelik@surrey.ac.uk}

\end{document}